\documentstyle[12pt,amssymb]{article} 

\font\csc=cmcsc10 scaled\magstep1
                
\font\teneufm=eufm10
\font\seveneufm=eufm7
\font\fiveeufm=eufm5
\newfam\eufmfam
\textfont\eufmfam=\teneufm
\scriptfont\eufmfam=\seveneufm
\scriptscriptfont\eufmfam=\fiveeufm

\let\goth\frak
\newcommand{\slth}{\widehat{\goth{sl}}_2}
\newcommand{\slt}{\goth{sl}_2}
\newcommand{\slnh}{\widehat{\goth{sl}}_n}

\newcommand{\g}{\goth{g}}
\newcommand{\cg}{\overline{\g}}
\newcommand{\uq}{U_q\bigl(\slth\bigr)}
\newcommand{\Aqp}{{\cal A}_{q,p}}

\newcommand{\Bqla}{{{\cal B}_{q,\lambda}}}
\newcommand{\Uqp}{U_{q,p}}

\font\fourteeneufm=eufm10 scaled\magstep2    
\font\seventeeneufm=eufm10 scaled\magstep3   
   
\newcommand{\slthbig}{\widehat{\mbox{\fourteeneufm sl}}_2}  
\newcommand{\slthBig}{\widehat{\mbox{\seventeeneufm sl}}_2} 
\newcommand{\gbig}{\mbox{\fourteeneufm g}} 
\newcommand{\gBig}{\mbox{\seventeeneufm g}} 

\newcommand{\Z}{{\Bbb Z}} 
\newcommand{\C}{{\Bbb C}} 

\newcommand{\F}{{\cal F}}
\newcommand{\hF}{{\hat{\F}}}

\newcommand{\cR}{{\cal R}}

\newcommand{\bR}{{\overline{R}}}
\newcommand{\hL}{\widehat{L}}

\newcommand{\tR}{{\widetilde{R}}}
\newcommand{\tL}{{\widetilde{L}}}
\newcommand{\ve}{\varepsilon}
\newcommand{\cL}{{\cal L}}

\newcommand{\la}{\lambda}
\newcommand{\La}{\Lambda}

\newcommand{\nn}{\nonumber}
\newcommand{\eqref}[1]{(\ref{#1})}
\newcommand{\bea}{\begin{eqnarray}}
\newcommand{\ena}{\end{eqnarray}}
\newcommand{\be}{\begin{eqnarray*}}
\newcommand{\en}{\end{eqnarray*}}
\newcommand{\lb}[1]{\label{#1}}

\newcommand{\tfrac}[2]{\textstyle \frac{#1}{#2}}
\newcommand{\End}{{\rm End}}

\newcommand{\id}{\hbox{id}}
\newcommand{\Ad}{\mbox{Ad}}
\newcommand{\br}[1]{{\langle #1 \rangle}} 

\newcommand{\ps}{{p^*}}
\newcommand{\rs}{{r^*}}

\def\infq4p#1{{(#1;q^4,p)_\infty}}

\newcommand{\EXP}[1]{{\exp\biggl\{#1\biggr\}}}

\newcommand{\tPsi}{\widetilde{\Psi}}
\newcommand{\tPhi}{\widetilde{\Phi}}
\newcommand{\pPhi}{\Phi'}
\newcommand{\pPsi}{{\Psi^*}'}

\newcommand{\hPhi}{\widehat{\Phi}}
\newcommand{\thPhi}{\widehat{\Phi}'}
\newcommand{\hPsi}{\widehat{\Psi}}
\newcommand{\thPsis}{\widehat{\Psi}^{*'}}
\newcommand{\hLp}{\widehat{L}^+}
\newcommand{\thL}{\widehat{L}^{+'}}
\newcommand{\hPsis}{\widehat{\Psi}^*}
\newcommand{\al}{\alpha}

\newcommand{\Remark}{\medskip \noindent {\it Remark.}\quad}
\newcommand{\qed}{\hfill \fbox{}\medskip}
\newcommand{\proof}{\medskip\noindent{\it Proof.}\quad }

\newtheorem{thm}{Theorem}[section]
\newtheorem{prop}[thm]{Proposition}
\newtheorem{lem}[thm]{Lemma}
\newtheorem{dfn}[thm]{Definition}
\newcommand{\ignore}[1]{}

\begin{document}
\font\csc=cmcsc10 scaled\magstep1

{\baselineskip=14pt
 \rightline{
 \vbox{
       \hbox{January 1998}\hfill 
       \hbox{HU-IAS/K-6}\\
       \phantom{a}\hfill\hbox{DPSU-98-2}
}}}

\vskip 11mm
\begin{center}

\hbox{\large\bf 
Elliptic algebra $\Uqp(\slthbig)$ : Drinfeld currents and vertex operators  
} 

\vskip11mm

{\csc Michio Jimbo}$\,^{1}$,
{\csc Hitoshi Konno}$\,^{2}$,
{\csc Satoru Odake}$\,{}^3$ and
{\csc Jun'ichi Shiraishi}$\,{}^4$
\\ 
{\baselineskip=15pt
\it\vskip.35in 
\setcounter{footnote}{0}\renewcommand{\thefootnote}{\arabic{footnote}}
\footnote{e-mail address : jimbo@kusm.kyoto-u.ac.jp}
Division of Mathematics, Graduate School of Science,\\
Kyoto University, Kyoto 606-8502, Japan\\
\vskip.1in 
\footnote{e-mail address : konno@mis.hiroshima-u.ac.jp}
Department of Mathematics, Faculty of Integrated Arts and Sciences,\\
Hiroshima University, Higashi-Hiroshima 739-8521,  Japan\\
\vskip.1in 
\footnote{e-mail address : odake@azusa.shinshu-u.ac.jp}
Department of Physics, Faculty of Science \\
Shinshu University, Matsumoto 390-8621, Japan\\
\vskip.1in 
\footnote{e-mail address : shiraish@momo.issp.u-tokyo.ac.jp}
Institute for Solid State Physics, \\
University of Tokyo, Tokyo 106-0032, Japan \\
}
\end{center}

\vskip5mm

\begin{abstract}
We investigate the structure of the elliptic algebra $\Uqp(\slth)$ 
introduced earlier by one of the authors. 
Our construction is based on a new set of generating series 
in the quantum affine algebra $U_q(\slth)$, 
which are elliptic analogs of the Drinfeld currents. 
They enable us to identify $\Uqp(\slth)$ with
the tensor product of $U_q(\slth)$ and 
a Heisenberg algebra generated by $P,Q$ with $[Q,P]=1$. 
In terms of these currents, 
we construct an $L$ operator satisfying the dynamical $RLL$ relation 
in the presence of the central element $c$. 
The vertex operators of Lukyanov and Pugai arise as 
`intertwiners' of $\Uqp(\slth)$ for the level one representation, 
in the sense to be elaborated on in the text. 
We also present vertex operators with higher level/spin 
in the free field representation. 
\end{abstract}

math.QA/9802002
(to appear in Communications in  Mathematical Physics)

\newpage

\setcounter{footnote}{0}
\renewcommand{\thefootnote}{\arabic{footnote})}

\setcounter{section}{0}
\setcounter{equation}{0}
\section{Introduction} \lb{sec:1}

\subsection{Vertex operators in SOS models} 

The principle of infinite dimensional symmetry 
has seen an impressive success in conformal field theory (CFT).
With the aim of understanding non-critical lattice models 
in the same spirit, 
the method of algebraic analysis  \cite{DFJMN,JMMN,JM} 
has been developed.
In this approach, 
a central role is played by the notion of vertex operators (VO's). 
There are two  kinds of VO's with distinct physical significance: 
the type I VO, which describes the operation of adding one lattice site, 
and the type II VO, which  
plays the role of particle creation/annihilation operators. 
In the most typical example of the XXZ spin chain, 
these VO's have a clear mathematical meaning as 
intertwiners \cite{FR} of certain modules 
over the quantum affine algebra $\uq$. 

An important class of CFT is the minimal unitary series \cite{BPZ}. 
Their lattice counterpart are 
the solvable models of Andrews-Baxter-Forrester (ABF) \cite{ABF}. 
These are `solid-on-solid' (SOS, or `face') models 
whose Boltzmann weights are expressed by elliptic functions. 
Their Lie theoretic generalizations have also been studied extensively
\cite{JMO,DJKMO3,JMO3}. 
The vertex operator approach to the ABF models and their fusion hierarchy 
was formulated in \cite{JMOh} by a coset-type construction. 
In \cite{JMOh}, 
$\uq$ was used only as an auxiliary tool to define the VO's, 
and its role as a symmetry algebra was somewhat indirect. 
In \cite{LukPug2}, Lukyanov and Pugai 
constructed a free boson realization of type I VO's for the ABF models. 
(The formulas for type II VO's can be found in \cite{MW96}.)
They have shown further that these VO's commute 
with the action of the deformed Virasoro algebra (DVA) \cite{qVir}, 
making clear the parallelism with CFT. 
However, unlike the case of CFT, 
the VO's did not allow for direct interpretation as intertwiners, 
because DVA lacks a coproduct\footnote{The usual coproduct for 
the Virasoro algebra has no non-trivial deformation.}.
It has remained an open problem to understand the conceptual meaning
of VO's. 

In \cite{Konno}, one of the authors introduced 
an elliptic algebra $\Uqp\bigl(\slth\bigr)$ 
and proposed it as 
an algebra of screening currents of conjectural extended DVA 
associated with the fusion SOS models.
The aim of the present article is to continue the 
study of  $\Uqp\bigl(\slth\bigr)$, and to show
that it offers a characterization of the VO's for SOS models 
in close analogy with the XXZ model.

\subsection{Face type elliptic algebras} 

Through an attempt to understand integrable models 
based on elliptic Boltzmann weights, 
various versions of `elliptic quantum groups' 
\cite{FIJKMY,Fel95,EF,Fron,Fron1} have been introduced. 
According to Fr\o nsdal \cite{Fron,Fron1}, elliptic quantum groups are
nothing but quantum affine algebras $U_q(\g)$ equipped with a 
coproduct different from the original one. 
The resulting objects are quasi-Hopf algebras in the sense of 
Drinfeld \cite{QHA}.  
Throughout this paper, we restrict our attention to the elliptic algebra of 
face type associated with $\g=\slth$, 
denoted as $\Bqla(\slth)$ in \cite{JKOS1}. 

In Fr\o nsdal's approach, the quasi-Hopf structures are defined by twistors 
given as formal series in the deformation parameters. 
An explicit construction for the twistors 
was given in \cite{JKOS1}.
(A very similar construction was presented independently in \cite{ABRR}.)
The $L$-operators and the VO's for the elliptic algebra 
can be obtained by `dressing' those of $\uq$ with the twistor  
(up to some subtleties about the fractional powers which will be 
discussed shortly). 
{}From this point of view, 
the construction of the VO's in bosonic representations is reduced to the 
determination of the image of the twistors. 
However, the solution of this issue is not known to us at this moment. 

For the bosonic realization of quantum affine algebras, 
the best suited presentation is in terms of the Drinfeld currents. 
In this paper we aim at an alternative construction of 
$L$ operators and VO's based on an elliptic analog of Drinfeld currents. 
These operators satisfy the same relations as those 
derived from the quasi-Hopf approach \cite{JKOS1}. 
Though the precise relation is not known, 
we expect that these two methods give equivalent answers. 
Our construction is inspired by the work of Enriquez and Felder \cite{EF}, 
who introduced Drinfeld-type currents defined on an elliptic curve 
and constructed the twistor by a quantum factorization method. 
The algebra $\Uqp\bigl(\slth\bigr)$ in \cite{Konno} 
and $U_{\hbar}\g(\tau)$ in \cite{EF} 
are both central extensions of the same algebra, 
but there are significant differences. 
We shall discuss more about this in section \ref{sec:6.2}. 

\subsection{Outline of the results}

Let us describe the content of this paper.
Our starting point is to introduce 
a new set of currents of $\uq$ carrying a parameter $r$ 
(the elliptic modulus), obtained by modifying the usual Drinfeld currents.
We shall refer to them as `elliptic currents'. 
They satisfy commutation relations with coefficients 
written in infinite products. 
The latter are essentially the Jacobi theta functions 
but not quite so, since the elliptic currents, 
and hence these coefficients, 
comprise only integral powers in the Fourier mode expansions. 
In order to have relations written in theta functions alone, 
we need to supply fractional powers. 
For this purpose we introduce 
`by hand' a pair of generators $P,Q$ which commute with $\uq$ 
and satisfy $[Q,P]=1$. 
Adjoining $P,Q$ to the elliptic currents, we obtain `total currents' 
whose commutation relations coincide with the defining 
relations of the algebra $\Uqp(\slth)$ \cite{Konno}
(see \eqref{Ucom1}-\eqref{Ucom10}). 
In other words, we can identify $\Uqp(\slth)$ with 
the tensor product of $\uq$ and the Heisenberg algebra generated by $P,Q$. 
The algebra $\Bqla(\slth)$ mentioned above
is the subalgebra of $U_{q,p}(\slth)$ isomorphic to $U_q(\slth)$,
and is equipped with a coproduct defined via the twistor.
However, this coproduct does {\it not} seem to extend naturally
to the full algebra $U_{q,p}(\slth)$.
That is, $\Bqla(\slth)$ is a quasi-Hopf algebra while
(to our knowledge) $U_{q,p}(\slth)$ is not.
We emphasize that
the intertwining relations for VO's will be based on the quasi-Hopf
structure of the former.

A characteristic feature of the elliptic algebras is that, 
in the presence of the central element $c$, 
we are forced to deal with 
two different elliptic moduli $r$ and $\rs=r-c$ simultaneously \cite{FIJKMY}. 
{}From Fr\o nsdal's point of view, it is an effect of the quasi-Hopf twisting. 
The appearance of two different curves 
makes it difficult to apply the geometric method of \cite{EF}. 
Instead, we take a more pedestrian approach. 
Motivated by similar formulas in \cite{EF}, 
we introduce `half currents' as certain 
contour integrals of the total currents.
They have an advantage that the coefficients of the 
commutation relations can be written solely in terms of theta functions
(as opposed to the delta functions appearing in the relations 
for the total currents). 
We then borrow the idea of the Gau{\ss} decomposition \cite{DF}
to compose an $L$ operator out of the half currents, and show that 
it satisfies the expected (dynamical) $RLL$ relation
(Proposition \ref{prop:RLL1}, \ref{prop:RLL2}). 

The construction of the 
$L$ operator allows us to study the VO's in the bosonic representation.
Let us first consider $\Bqla(\slth)$.
As is clear from the construction, 
the elliptic currents can be realized in the same 
bosonic Fock spaces as with the Drinfeld currents of $\uq$. 
(We regard $p=q^{2r}$ as a formal parameter.) 
The VO's for $\Bqla(\slth)$ are a family of intertwiners 
$\Phi(z,s)$, $\Psi^*(z,s)$ carrying a parameter $s$, and 
their intertwining relations involve a shift of $s$ 
(see \eqref{int1}-\eqref{int2}). 
With the adjunction of $P,Q$, 
the algebra $\Uqp(\slth)$ has an enlarged Fock module. 
It has the decomposition $\F=\bigoplus_s \F_s$ into eigenspaces 
of $P$, each eigenspace $\F_s$ being a Fock module for $\uq$. 
Accordingly we modify further the VO's with $Q$, 
\bea
&&\hPhi(u,s)=z^{\frac{1}{2r}(\frac{1}{2}{h^{(2)}}^2+(s+h^{(1)})h^{(2)})}
\Phi(q^cz,s),
\lb{VOU1}\\
&&\hPsi^*(u,s)=
\Psi^*(z,s)z^{-\frac{1}{2\rs}(\frac{1}{2}{h^{(1)}}^2+sh^{(1)})}
e^{Qh^{(1)}},
\lb{VOU2}
\ena
where $z=q^{2u}$, $h^{(1)}=h\otimes 1, h^{(2)}=1\otimes h$, $h$ being 
the `Cartan' generator of $\uq$. 
Solving the intertwining relations for level one, 
we find that the VO's of Lukyanov and Pugai arise 
in the form \eqref{VOU1}, \eqref{VOU2}, 
apart from certain signs in the intertwining relations for $\Phi(z,s)$ 
and $\Psi^*(z,s)$. 
(For the discussion of the signs, see subsections 
\ref{subsec:5.2} and \ref{subsec:5.4}.) 
We also calculate formulas for VO's associated with 
higher spin representations. 

\subsection{Plan of the text} 

The text is organized as follows. 

In section 2, we recall some results of \cite{JKOS1} 
which are relevant to the following sections. 
In Section3, we introduce the elliptic currents of $\uq$, 
and discuss its relation to $\Uqp(\slth)$. 
In section 4, we introduce the `half currents' and derive their 
commutation relations. 
We then arrange them in the form of a Gau{\ss} decomposition 
to define the $L$ operator. 
In section 5, we describe the VO's in the bosonic representation. 
In Section 6 we discuss the connection to other works and mention 
some open problems. 
The text is followed by four appendices. 
In appendix \ref{app:1}, we give the elliptic currents 
for the general non-twisted affine Lie algebra $\g$. 
In Appendix \ref{app:2} we discuss an 
elliptic analog of the Drinfeld coproduct, 
and show that it also arises as 
a quasi-Hopf twist from the usual Drinfeld coproduct.
In Appendix \ref{app:3} we study the 
evaluation modules and $R$ matrix in the spin $l/2$ representation. 
Finally, in Appendix \ref{app:4} we review the free field realization of 
the algebra $\Uqp(\slth)$. 

While preparing this manuscript, we became aware 
of the paper by Hou et al. \cite{HouY} 
which has some overlap with the content of the present paper. 


\setcounter{section}{1}
\setcounter{equation}{0}
\section{$RLL$ and intertwining relations}\lb{sec:2}

The purpose of this section is to 
set up the form of the $RLL$ and intertwining relations 
which we are going to study. 

\subsection{Previous results}\lb{sec:2.1}

In order to fix the notation, let us recall the results of \cite{JKOS1}
relevant to the present paper. 
We consider the quantum affine algebra $U_q=\uq$ with standard generators 
$e_i,f_i,h_i$ ($i=0,1$) and $d$. 
The canonical central element is $c=h_0+h_1$. 
We retain the convention of \cite{JKOS1} for the coproduct $\Delta$, 
though the details are not necessary here. 

Henceforth we shall write $h=h_1$. 
In \cite{JKOS1}, we have constructed a twistor $F(\la)\in U_q^{\otimes 2}$. 
Changing slightly the notation, 
let us set $\la=(\rs+2)d+(s+1)\frac{1}{2}h$ and write $F(\la)$ as $F(\rs,s)$. 
Then $F(\rs,s)$ is a formal power series in $q^{2(\rs-s)}$ and $q^{2s}$, 
satisfying the shifted cocycle condition
\bea
F^{(12)}(\rs,s)\left(\Delta\otimes\id\right)F(\rs,s)
=
F^{(23)}(\rs+c^{(1)},s+h^{(1)})\left(\id\otimes\Delta\right)F(\rs,s).
\lb{cocy}
\ena
We obtain the quasi-Hopf algebra $\Bqla=\Bqla(\slth)$
by twisting $U_q$ via this $F$.
Here and after, the superscripts refer to the tensor components; 
for instance, $F^{(23)}=1\otimes F$, $h^{(1)}=h\otimes 1\otimes 1$. 
Let $\cR$ be the universal $R$ matrix of $U_q$.
The `dressed' $R$ matrix 
$\cR(\rs,s)=F^{(21)}(\rs,s)\cR F^{(12)}(\rs,s)^{-1}$ 
of $\Bqla$ satisfies the dynamical YBE, 
\bea
&&\cR^{(12)}(\rs+c^{(3)},s+h^{(3)})\cR^{(13)}(\rs,s)
\cR^{(23)}(\rs+c^{(1)},s+h^{(1)})
\nn\\
&&\quad=
\cR^{(23)}(\rs,s)\cR^{(13)}(\rs+c^{(2)},s+h^{(2)})\cR^{(12)}(\rs,s).
\lb{DYBE}
\ena

Let $(\pi_V,V)$ be a finite dimensional $U'_q$-module 
where $U_q'$ is the subalgebra generated by $e_i,f_i,h_i$ ($i=0,1$). 
Let $(\pi_{V,z},V_z)$ denote the evaluation module
\be
&&\pi_{V,z}(a)=\pi_V\circ\Ad(z^d)(a) \quad(a\in U_q'), 
\qquad
\pi_{V,z}(d)=z\frac{d}{dz},
\qquad 
V_z=V[z,z^{-1}].
\en
Setting 
\bea
&&R^+_{VW}(z_1/z_2,\rs,s)=\left(\pi_{V,z_1}\otimes\pi_{W,z_2}\right)
\cR(\rs,s),
\lb{Rp}\\
&&
L^+_{V}(z,\rs,s)=\left(\pi_{V,z}\otimes\id\right)
q^{c\otimes d+d\otimes c}\cR(\rs,s),
\lb{Lp}
\ena
we have the dynamical $RLL$ relation for $\Bqla$,
\bea
&&R^{+(12)}_{VW}(z_1/z_2,\rs+c,s+h)
L^{+(1)}_{V}(z_1,\rs,s)L^{+(2)}_{W}(z_2,\rs,s+h^{(1)})
\nn\\
&&\quad =
L^{+(2)}_{W}(z_2,\rs,s)L^{+(1)}_{V}(z_1,\rs,s+h^{(2)})
R^{+(12)}_{VW}(z_1/z_2,\rs,s). 
\lb{RLL}
\ena
Hereafter we shall write 
\bea
&&r=\rs+c,
\qquad 
R^+_{VW}(z,s)=R^{+}_{VW}(z,r,s), 
\quad
R^{*+}_{VW}(z,s)=R^{+}_{VW}(z,\rs,s),
\lb{RRs}
\ena
and normally suppress the $\rs$-dependence.
In this paper we will not consider the $L^-$ operator since 
it can be obtained from $L^+$, 
see Proposition 4.3 in \cite{JKOS1}.

Let now $\F,\F'$ be highest weight $U_q$-modules on which $c$ acts as 
a scalar $k$. 
Suppose we have intertwiners of $U_q$-modules 
\be
&&\Phi_V(z): \F \longrightarrow \F'\otimes V_z,
\\
&&\Psi^*_V(z):V_z\otimes\F \longrightarrow \F', 
\en
which we refer to as vertex operators (VO's) of type I and type II, 
respectively. 
Then the `dressed' VO's for $\Bqla$
\bea
&&\Phi_V(z,s)=\left(\id\otimes \pi_{V,z}\right)F(\rs,s)\circ \Phi_V(z),
\lb{Ph}\\
&&
\Psi^*_V(z,s)=\Psi^*_V(z)\circ \left(\pi_{V,z}\otimes\id\right)F(\rs,s)^{-1}
\lb{Ps}
\ena
satisfy the following intertwining relations with the $L$ operators:
\bea
&&\Phi_{W}(q^kz_2,s)L^+_{V}(z_1,s)
=R^+_{VW}(z_1/z_2,s+h)L^+_{V}(z_1,s)
\Phi_{W}(q^kz_2,s+h^{(1)}),
\lb{int1}\\
&&
L^+_{V}(z_1,s)\Psi^*_{W}(z_2,s+h^{(1)})
=\Psi^*_{W}(z_2,s)L^+_{V}(z_1,s+h^{(2)})
R^{*+}_{VW}(z_1/z_2,s).
\lb{int2}
\ena
We note that all the operators \eqref{Rp}, \eqref{Lp},
\eqref{Ph}, \eqref{Ps} are 
formal Laurent series comprising only integral powers of $z$. 

\Remark 
In the present paper, we shall adopt the universal $R$ matrix 
$\cR={\cR^{'(21)}}^{-1}$, 
where $\cR'$ is given in (2.8) of \cite{JKOS1}. 
This is purely a matter of convention. 
The properties (2.11)-(2.14) holds equally well for 
$\cR$ and $\cR'$, and hence the same construction applies. 

\subsection{Fractional powers}\lb{sec:2.2}

The $RLL$-relation \eqref{RLL} is unchanged under the transformation 
of the form 
\bea
&&L^{+'}_V(z,s)=\mu_V(s,h^{(1)})\mu(s+h^{(1)},h)L^+_V(z,s)
\mu^{*}_V(s+h^{(2)},h^{(1)})^{-1}\mu(s,h)^{-1},
\lb{Lpr}\\
&&
R^{'+}_{VW}(z_1/z_2,s)=\mu_V(s,h^{(1)})\mu_W(s+h^{(1)},h^{(2)})
\nn\\
&&\qquad\qquad\qquad\qquad \times 
R^+_{VW}(z_1/z_2,s)
\mu_V(s+h^{(2)},h^{(1)})\mu_W(s,h^{(2)}).
\lb{Rpr}
\ena
Here $\mu_V(s,h)$, $\mu_W(s,h)$, $\mu(s,h)$ are functions 
possibly depending on $z$ and $r$, and $\mu^*_V(s,h)$ means 
$\mu_V(s,h)\bigl|_{r\rightarrow \rs}$. 
This corresponds to the freedom of 
changing the twistor by `shifted coboundary'. 
Exploiting this freedom, we 
modify the $R$ matrix by a fractional power of $z$
so that it can be expressed in terms of the Jacobi theta functions.

Consider the image of the $R$ matrix in the 
evaluation modules $(\pi_{l,z},V_{l,z})$ of spin $l/2$ 
(see Appendix \ref{app:3})
\be
R^+_{lm}(z_1/z_2,s)
=\left(\pi_{l,z_1}\otimes\pi_{m,z_2}\right)\cR(r,s).
\en
In Appendix \ref{app:3}, we show that it has the form 
$R^+_{lm}(z,s)=\rho^+_{lm}(z,p)\bR_{lm}(z,s)$,
where $\rho^+_{lm}(z,p)$ is a scalar factor 
given in \eqref{rhop} and 
$\bR_{lm}(z,s)$ has the transformation property \eqref{pz}.
Write $z=q^{2u}$ and set 
\bea
&&\tR^+_{lm}(u,s)
=z^{\frac{1}{2r}(\frac{1}{2}{h^{(1)}}^2+(s+h^{(2)})h^{(1)})}
R^+_{lm}(z,s)
z^{-\frac{1}{2r}(\frac{1}{2}{h^{(1)}}^2+sh^{(1)})}.
\lb{tR}
\ena
The $R$ matrix \eqref{tR} comprises fractional powers of $z$, 
but (up to a scalar factor) becomes completely periodic, 
\be
\frac{1}{\rho^+_{lm}(pz,p)}\tR^+_{lm}(u+r,s)
=
\frac{1}{\rho^+_{lm}(z,p)}\tR^+_{lm}(u,s).
\en
It turns out that the $R$ matrix \eqref{tR} 
is expressible in terms of the Jacobi theta functions. 
An explicit expression for the case $l=m=1$ will be given 
in \eqref{Rmat} below. 

In accordance with \eqref{tR}, we modify the $L$ operator and 
VO's as
\bea
&&\tL^+_m(u,s)=z^{\frac{1}{2r}(\frac{1}{2}{h^{(1)}}^2+(s+h)h^{(1)})}
L^+_{V_m}(z,s)z^{-\frac{1}{2\rs}(\frac{1}{2}{h^{(1)}}^2+sh^{(1)})},
\lb{tL}\\
&&
\tPhi_l(u,s)=
z^{\frac{1}{2r}(\frac{1}{2}{h^{(2)}}^2+(s+h^{(1)})h^{(2)})}
\Phi_{V_l}(q^cz,s),
\lb{tPhi}\\
&&\tPsi^*_{l}(u,s)=\Psi^*_{V_l}(z,s)
z^{-\frac{1}{2\rs}(\frac{1}{2}{h^{(1)}}^2+sh^{(1)})}.
\lb{tPsi}
\ena
We shall focus attention to the $L$ operators $\tL^+(u,s)=\tL^+_{1}(u,s)$ 
associated with the spin $1/2$ representation. 
The following are consequences of \eqref{RLL}, \eqref{int1}, \eqref{int2}:
\bea
&&\tR^{+(12)}_{11}(u_1-u_2,s+h)
\tL^{+(1)}(u_1,s)\tL^{+(2)}(u_2,s+h^{(1)})
\nn\\
&&\quad =
\tL^{+(2)}(u_2,s)\tL^{+(1)}(u_1,s+h^{(2)})
\tR^{*+(12)}_{11}(u_1-u_2,s),
\lb{RLL2}\\
&&
\tPhi_{l}(u_2,s)\tL^+(u_1,s)
=\tR^+_{1l}(u_1-u_2,s+h)\tL^+(u_1,s)\tPhi_{l}(u_2,s+h^{(1)} ),
\lb{int3}\\
&&
\tL^+(u_1,s)\tPsi^*_{l}(u_2,s+h^{(1)})
=\tPsi^*_{l}(u_2,s)\tL^+(u_1,s+h^{(2)})\tR^{*+}_{1l}(u_1-u_2,s).
\lb{int4}
\ena
We shall study 
the relations \eqref{RLL2}--\eqref{int4} in the following sections. 


\setcounter{section}{2}
\setcounter{equation}{0}
\section{Elliptic currents and $\Uqp\bigl(\slthBig\bigr)$} \lb{sec:3}

In this section, we introduce Drinfeld-type currents of 
$U_q=\uq$ satisfying `elliptic' commutation relations. 
We then relate them to the elliptic algebra 
$\Uqp\bigl(\slth\bigr)$ of \cite{Konno} 
by adjoining a pair of generators $P,Q$ with $[Q,P]=1$. 

\subsection{Elliptic currents of $U_q(\slthbig)$}\lb{subsec:3.1}

First let us recall the Drinfeld currents of $U_q$ \cite{Dri88}. 
Hereafter we fix a complex number $q\neq 0$, $|q|<1$. 
We use the standard symbols 
\be
&&[n]=\frac{q^n-q^{-n}}{q-q^{-1}}.
\en
Let 
 $x^\pm_n$ ($n\in\Z$),  $a_n$ ($n\in\Z_{ \neq 0}$), $h$, $c$, $d$ 
denote the Drinfeld generators of  $U_q$. 
In terms of the generating functions 
\bea
&&x^\pm(z)=\sum_{n\in \Z}x^\pm_n z^{-n},
\lb{Dcur1}\\
&&\psi(q^{c/2}z)=q^h
\exp\left( (q-q^{-1}) \sum_{n>0} a_{n}z^{- n}\right),
\lb{Dcur2}\\
&&
\varphi(q^{-c/2}z)=q^{-h}
\exp\left(-(q-q^{-1})\sum_{n>0} a_{-n}z^{ n}\right), 
\lb{Dcur3}
\ena
the defining relations read as follows: 
\be
&&c :\hbox{ central },\\
&& [h,d]=0,\quad [d,a_n]=n a_n,\quad [d,x^{\pm}_n]=n x^{\pm}_n, \\
&&[h,a_n]=0,\qquad [h, x^\pm(z)]=\pm 2 x^{\pm}(z),\\
&&
[a_n,a_m]=\frac{[2n][c n]}{n}
q^{-c|n|}\delta_{n+m,0},\\
&&
[a_n,x^+(z)]=\frac{[2n]}{n}q^{-c|n|}z^n x^+(z),\\
&&
[a_n,x^-(z)]=-\frac{[2n]}{n} z^n x^-(z),
\\
&&(z-q^{\pm 2}w)
x^\pm(z)x^\pm(w)= (q^{\pm 2}z-w) x^\pm(w)x^\pm(z),
\\
&&[x^+(z),x^-(w)]=\frac{1}{q-q^{-1}}
\left(\delta\bigl(q^{-c}\frac{z}{w}\bigr)\psi(q^{c/2}w)
-\delta\bigl(q^{c}\frac{z}{w}\bigr)\varphi(q^{-c/2}w)
\right).
\en

We now introduce a new parameter $p$ and 
 modify \eqref{Dcur1}-\eqref{Dcur3} to define another set of currents. 
For notational convenience,  we will frequently write 
\be
&&p=q^{2 r}, \qquad 
p^*=p q^{-2c}=q^{2\rs}\qquad (\rs=r-c).
\en
Let us introduce two currents $u^\pm(z,p)\in U_q$ depending on $p$ by 
\bea
&&
u^+(z,p)=\exp\left(\sum_{n>0}\frac{1}{[r^*n]}a_{-n}(q^rz)^n\right),
\lb{dressup}\\
&&
u^-(z,p)=\exp\left(-\sum_{n>0}\frac{1}{[rn]}a_{n}(q^{-r}z)^{-n}\right).
\lb{dressum}
\ena

\begin{dfn}[Elliptic currents]
We define the currents $e(z,p)$, $f(z,p)$, $\psi^\pm(z,p)$  by
\bea
&& e(z,p)=u^+(z,p)x^+(z), 
\lb{dress1}
\\
&&f(z,p)=x^-(z)u^-(z,p),
\lb{dress2}
\\
&&\psi^+(z,p)
=u^+(q^{c/2} z,p)\psi(z)u^-(q^{-c/2}z,p),
\lb{dress3}
\\
&&\psi^{-}(z,p)
=u^+(q^{-c/2}z,p)\varphi(z)u^-(q^{c/2}z,p). 
\lb{dress4}
\ena
\end{dfn}
We will often drop $p$, and write $e(z,p)$ as $e(z)$ and so forth. 

The merit of these currents is that they obey the following `elliptic' 
commutation relations.
\begin{prop} \label{ellipticcom}
\bea
&&
\psi^{\pm}(z)\psi^{\pm}(w)
=\frac{\Theta_{p}\left(q^{-2}z/w\right)}{\Theta_{p}\left(q^2 z/w\right)}
\frac{\Theta_{p^*}\left(q^{2}z/w\right)}{\Theta_{p^*}\left(q^{-2} z/w\right)}
\psi^{\pm}(w)\psi^{\pm}(z),
\lb{Drcom1}\\
&&\psi^{+}(z)\psi^{-}(w)
=\frac{\Theta_{p}\left(pq^{-c-2}z/w\right)}
{\Theta_{p}\left(pq^{-c+2} z/w\right)}
\frac{\Theta_{p^*}\left(p^*q^{c+2}z/w\right)}
{\Theta_{p^*}\left(p^*q^{c-2} z/w\right)}
\psi^{-}(w)\psi^{+}(z),
\lb{Drcom2}\\
&&\psi^{\pm}(z)e(w)\psi^{\pm}(z)^{-1}=
q^{-2}\frac{\Theta_{p^*}(q^{\pm c/2+2} z/w)}{\Theta_{p^*}(q^{\pm c/2-2} z/w)}
e(w),
\lb{Drcom3}\\
&&\psi^{\pm}(z)f(w)\psi^{\pm}(z)^{-1}=
q^2\frac{\Theta_{p}(q^{\mp c/2-2} z/w)}{\Theta_{p}(q^{\mp c/2+2} z/w)}
f(w),
\lb{Drcom4}\\
&&[e(z),f(w)]=\frac{1}{q-q^{-1}}
\left(\delta\bigl(q^{-c}\frac{z}{w}\bigr)\psi^+(q^{c/2}w)
-\delta\bigl(q^{c}\frac{z}{w}\bigr)\psi^-(q^{-c/2}w)
\right).
\lb{Drcom5}
\ena
\end{prop}
Here we have used the standard symbols 
\be
&&\Theta_p(z)=(z;p)_\infty(pz^{-1};p)_\infty(p;p)_\infty,
\\
&&(z;t_1,\cdots,t_k)_\infty=
\prod_{n_1,\cdots,n_k\ge 0}(1-zt_1^{n_1}\cdots t_k^{n_k}).
\en

It will become convenient later to consider also the current 
\bea
&&
k(z)=\exp\left(\sum_{n>0}\frac{[n]}{[2n][r^*n]}a_{-n} (q^c z)^n\right)
\exp\left(-\sum_{n>0}\frac{[n]}{[2n][rn]}a_n z^{-n}\right). 
\lb{dress5}
\ena
The $\psi^\pm(z)$ are related to $k(z)$ by the formula 
\bea
&& \psi^\pm(p^{\mp(r-\frac{c}{2})}z)=\kappa q^{\pm h}k(q z)k(q^{-1}z),  
\label{eq:psipm}
\\
&&\kappa=\frac{\xi(z;p^*,q)}{\xi(z;p,q)} \Biggl|_{z=q^{-2}},
\label{kappac}
\ena
where the function 
\bea
\xi(z;p,q)=
\frac{(q^2z;p,q^4)_\infty(pq^2z;p,q^4)_\infty}
{(q^4z;p,q^4)_\infty(p z;p,q^4)_\infty}
\lb{xi}
\ena
is a solution of the difference equation 
\be
\xi(z;p,q)\xi(q^2z;p,q)=\frac{(q^2z;p)_\infty}{(pz;p)_\infty}.
\en
We have the commutation relations supplementing \eqref{Drcom1}-\eqref{Drcom5},
\bea
&&
k(z)k(w)=\frac{\xi(w/z;p,q)}{\xi(w/z;p^*,q)}
\frac{\xi(z/w;p^*,q)}{\xi(z/w;p,q)}k(w)k(z),
\lb{Drcom6}\\
&&k(z)e(w)k(z)^{-1}
=\frac{\Theta_{p^*}\left(p^{*1/2}qz/w\right)}
{\Theta_{p^*}\left(p^{*1/2}q^{-1}z/w\right)} e(w),
\lb{Drcom7}\\
&&k(z)f(w)k(z)^{-1}
=\frac{\Theta_{p}\left(p^{1/2}q^{-1}z/w\right)}
{\Theta_{p}\left(p^{1/2}q z/w\right)}f(w).
\lb{Drcom8}
\ena

The commutation relations \eqref{Drcom1}-\eqref{Drcom5} have been 
proposed earlier in \cite{DI97}, but the direct connection 
with the usual Drinfeld currents was not known. 
In Appendix \ref{app:2} we discuss also 
the Drinfeld type coproduct 
for the elliptic currents \eqref{dress1}-\eqref{dress4}, \eqref{dress5}. 

\Remark Strictly speaking, the currents \eqref{dressup}-\eqref{dress4},
\eqref{dress5} are generating series whose coefficients belong to 
a completion of $\uq\otimes \C[[p]]$.
At this level, $p$ should be treated as an indeterminate. 
However, in the concrete representations we are going to discuss, 
such as evaluation modules and Fock modules, these currents have also 
analytical meaning. 
(For a formula for the currents in spin $l/2$ evaluation modules, 
see Appendix \ref{app:3}.)
We will not go into this point any further, and later treat 
$p$, $p^*$ as complex numbers satisfying $|p|,|p^*|<1$. 
\medskip

\subsection{Elliptic algebra $\Uqp\bigl(\slthbig\bigr)$}\label{subsec:3.2}

The elliptic algebra $\Uqp\bigl(\slth\bigr)$ in \cite{Konno} 
is very similar to the algebra of the elliptic currents 
\eqref{Drcom1}-\eqref{Drcom5}, \eqref{Drcom6}-\eqref{Drcom8}. 
In the former, the coefficients of the relations 
are written in terms of the Jacobi elliptic theta function, 
which differs from $\Theta_p(z)$ used in the latter by a 
simple factor (see \eqref{theta} below).
Let us discuss the precise connection between the two algebras. 

For this purpose, it is more convenient to work with the `additive' notation. 
Following \cite{Konno}, we use the parameterization 
\be
&&q=e^{-\pi i/r\tau},
\\
&&p=e^{-2\pi i/\tau}, \quad
p^*=e^{-2\pi i/\tau^*}
\qquad (r\tau=r^*\tau^*),
\\
&&z=q^{2u}=e^{-2\pi i u/r\tau}. 
\en
We also use the Jacobi theta functions
\bea
\theta(u)=
q^{{u^2 \over r}-u}
{ \Theta_p(q^{2 u}) \over (p;p)_\infty^3},
\quad
\theta^*(u)=
q^{{u^2 \over r^*}-u}
 {\Theta_{p^*}(q^{2 u}) \over (p^*;p^*)_\infty^3}.
\lb{theta}
\ena
The function $\theta(u)$ has a zero at $u=0$, 
enjoys the quasi-periodicity property 
\be
\theta(u+r)=-\theta(u),
\qquad 
\theta(u+r\tau)=-e^{-\pi i \tau-\tfrac{2\pi i u}{r}}\theta(u), 
\en
and is so normalized that 
\begin{eqnarray*}
\oint_{C_0} {dz \over 2 \pi i z} {1 \over \theta(-u)}=1,
\end{eqnarray*}
where $C_0$ is a simple closed curve in the $u$-plane 
encircling $u=0$ counterclockwise. 
The same holds for $\theta^*(u)$, 
with $r$ and $\tau$ replaced by $\rs$ and $\tau^*$ respectively. 

Now let us introduce new generators $P,Q$ such that 
\bea
[Q,P]=1,
\qquad
\hbox{ $Q,P$ commute with $\uq$}.
\lb{QP}
\ena
With the aid of them, we define the `total' currents obtained by modifying 
the elliptic currents of the previous subsection \ref{subsec:3.1}.  
Below we shall use the notation for the conformal weight 
\bea
\Delta_{l,r}={l(l+2) \over 4r}.
\lb{confwt}
\ena

\begin{dfn}[Total currents]
We define the currents $K(u),E(u),F(u),H^{\pm}(u)$ by 
\bea
&&K(u)=k(z) e^Q  z^{\Delta_{-P-h-1}-\Delta_{-P-h}-
          \Delta^*_{-P-1}+\Delta^*_{-P}},
\lb{U1}
\\
&&E(u)=e(z) e^{2 Q}z^{-\Delta^*_{-P-1}+\Delta^*_{-P+1}},
\lb{U2}\\
&&F(u)= f(z) z^{\Delta_{-P-h-1}-\Delta_{-P-h+1}},
\lb{U3}
\\
&&H^{\pm}(u)
=\psi^\pm(z)e^{2Q}
\left(q^{\pm \bar{r}}z\right)^{\Delta_{-P-h-1}-\Delta_{-P-h+1}
-\Delta^*_{-P-1}+\Delta^*_{-P+1}}.
\lb{U4}
\ena
Here we have set $z=q^{2u}$, 
$\Delta_l=\Delta_{l,r}$, $\Delta^*_l=\Delta_{l,\rs}$, 
and $\bar{r}=r-{c \over 2}$. 
\end{dfn}
The currents $K(u)$ and $H^\pm(u)$ are related by 
\bea
&&H^{\pm}(u)=\kappa K\left(u\pm \frac{\bar{r}}{2}+\frac{1}{2}\right)
K\left(u\pm \frac{\bar{r}}{2}-\frac{1}{2}\right),
\\
&&H^-(u)=H^+(u-\bar{r}), 
\ena
with the same $\kappa$ as in \eqref{kappac}.
We shall refer to \eqref{U1}-\eqref{U4} as total currents. 

{}From the commutation relations of the elliptic currents, 
we can derive those of the total currents. 
Let us introduce a function $\rho(u)$ by 
\bea
&&\rho(u)=\frac{\rho^{+*}(u)}{\rho^+(u)},
\lb{rho1}
\ena
where
\bea
&&\rho^+(u)=z^{\frac{1}{2r}}\rho^+_{11}(z,p)
=z^{\frac{1}{2r}}q^{\frac12}
\frac{\{pq^2z\}^2}{\{pz\}\{pq^4z\}}
\frac{\{z^{-1}\}\{q^4z^{-1}\}}{\{q^2z^{-1}\}^2},
\lb{rho11}\\
&&\{z\}=(z;p,q^4)_\infty,
\lb{pq4}
\ena
$\rho^+_{lm}(z,p)$ is given in \eqref{rhop}, 
and $\rho^{+*}(u)=\rho^+(u)\bigl|_{r\rightarrow \rs}$.

\begin{prop}
The following commutation relations hold:
\bea
\!\!\!\!\!\!\!\!\!\!&&K(u)K(v)=\rho(u-v)K(v)K(u),
\lb{Ucom1}\\
\!\!\!\!\!\!\!\!\!\!&&K(u)E(v)=\frac{\theta^*(u-v+\frac{1-r^*}{2})}
{\theta^*(u-v-\frac{1+r^*}{2})}E(v)K(u),
\lb{Ucom2}\\
\!\!\!\!\!\!\!\!\!\!&&K(u)F(v)=\frac{\theta(u-v-\frac{1+r}{2})}
{\theta(u-v+\frac{1-r}{2})}F(v)K(u),
\lb{Ucom3}\\
\!\!\!\!\!\!\!\!\!\!&&E(u)E(v)=\frac{\theta^*(u-v+1)}{\theta^*(u-v-1)}E(v)E(u),
\lb{Ucom4}\\
\!\!\!\!\!\!\!\!\!\!&&F(u)F(v)=\frac{\theta(u-v-1)}{\theta(u-v+1)}F(v)F(u),
\lb{Ucom5}\\
\!\!\!\!\!\!\!\!\!\!&&[E(u),F(v)]=\frac{1}{q-q^{-1}}
\left(\delta(u-v-\frac{c}{2})H^+(u-\frac{c}{4})
-\delta(u-v+\frac{c}{2})H^-(v-\frac{c}{4})\right),
\lb{Ucom6}
\\
\!\!\!\!\!\!\!\!\!\!&&H^{+}(u)H^{-}(v)=\frac{\theta(u-v-\frac{c}{2}-1)}
{\theta(u-v-\frac{c}{2}+1)}
\frac{\theta^*(u-v+\frac{c}{2}+1)}
{\theta^*(u-v+\frac{c}{2}-1)}H^{-}(v)H^{+}(u),
\lb{Ucom7}\\
\!\!\!\!\!\!\!\!\!\!&&H^{\pm}(u)H^{\pm}(v)=\frac{\theta(u-v-1)}{\theta(u-v+1)}
\frac{\theta^*(u-v+1)}{\theta^*(u-v-1)}H^{\pm}(v)H^{\pm}(u),
\lb{Ucom8}\\
\!\!\!\!\!\!\!\!\!\!&&H^{\pm}(u)E(v)=
\frac{\theta^*(u-v\pm\frac{c}{4}+1)}
{\theta^*(u-v\pm\frac{c}{4}-1)}E(v)H^{\pm}(u),
\lb{Ucom9}\\
\!\!\!\!\!\!\!\!\!\!&&H^{\pm}(u)F(v)=
\frac{\theta(u-v\mp\frac{c}{4}-1)}
{\theta(u-v\mp\frac{c}{4}+1)}F(v)H^{\pm}(u). 
\lb{Ucom10}
\ena
Here $\delta(u)$ means 
$\displaystyle \delta(z)=\sum_{n\in\Z} z^n$ ($z=q^{2u}$).
\end{prop}
It is in this form that the algebra $\Uqp\bigl(\slth\bigr)$ was presented 
in \cite{Konno}. 
Thus we arrive at the following interpretation: 
\begin{dfn}[Algebra $\Uqp\bigl(\slth\bigr)$]
We define the algebra $\Uqp\bigl(\slth\bigr)$ to be the
tensor product of $\uq$ and a Heisenberg algebra with generators
$Q,P$ \eqref{QP}.
\end{dfn}
We note that, when $c=0$, the commutation 
relations \eqref{Ucom1}-\eqref{Ucom10} 
also coincide with those of Enriquez-Felder \cite{EF} with $K=0$. 

In the bosonization of Appendix D, 
the elements $P$, $h$ and $Q$ are given as follows ($c=k$):
\be
&&P-1=\widehat{\Pi}=\sqrt{\frac{2r(r-k)}{k}}P_0-\frac{r-k}{k}P_2,
\\
&&P+h-1=\widehat{\Pi}'=\sqrt{\frac{2r(r-k)}{k}}P_0-\frac{r}{k}P_2,
\\
&&Q=-\sqrt{\frac{k}{2r(r-k)}}iQ_0 =-\sqrt{2}\alpha_0iQ_0,
\en
where $\widehat{\Pi},\widehat{\Pi}'$ are the notation of \cite{Konno}.
The physical meaning of the quantity  $2\alpha_0$ is 
the anomalous charge of the boson field $\phi_0$.
We note also that, using the notation \eqref{confwt} of the conformal weight, 
the $L_0$ operator in \cite{Konno} can be written as 
\bea
&&L_0=-d+\Delta_{-P+1,r-k}-\Delta_{-P-h+1,r}. 
\lb{Lzero}
\ena

The elliptic currents of $U_q(\slth)$ and the algebra $U_{q,p}(\slth)$
are naturally extended to those associated with arbitrary non-twisted 
affine Lie algebras. In Appendix A, we give a summary of the results
and discuss their significance.


\setcounter{section}{3}
\setcounter{equation}{0}
\section{The $RLL$ relations}\lb{sec:4}

One of the goals of the present paper is to describe the
vertex operators (VO's) $\Phi_V(z,s)$, $\Psi^*_V(z,s)$ in
the bosonic representation of the algebras $\Bqla(\slth)$ 
and $U_{q,p}(\slth)$ given in Appendix D.
The intertwining relations for the VO's \eqref{int1},\eqref{int2} 
are based on the operator $L^+_V(z,r^*,s)$ defined in \eqref{Lp}.
In order to compute the VO's, therefore, we need the image of the 
`dressed' universal $R$ matrix ${\cal R}(r^*,s)$ in the Fock space.
The latter is given as an infinite product of the universal $R$
matrix for $U_q(\slth)$ \cite{JKOS1},
but we do not know how to calculate it at this moment.

In this section, we take an alternative approach.
Namely we utilize the elliptic currents to construct a $2\times 2$ 
matrix operator $L^+(u,P)$ (see \eqref{Gauss}, \eqref{Lplus}),
and show that it satisfies the same $RLL$-relation \eqref{RLL}
as for $L^+_V(z,r^*,s)$ with $V$ being the spin $1/2$ representation.
Though we do not know a proof, from this construction
we expect that (modulo perhaps some base change) this $L^+(u,P)$ is
the same as $L^+_V(z,r^*,s)$ of $\Bqla(\slth)$ (with $s=P$).

\subsection{Half currents}\lb{subsec:4.1}

The commutation relations of the total currents $E(u)$ and $F(u)$ 
involve delta functions. 
We are going to modify them so as to have commutation relations 
involving only `ordinary' functions. 
Motivated by a similar construction in \cite{EF}, 
we define the half currents of $\Uqp\bigl(\slth\bigr)$ as follows. 

\begin{dfn}[Half currents]
We set 
\bea
&&
K^+(u)=K(u+\tfrac{r+1}{2}),
\label{kplush}\\
&&E^+(u)
=a^* \oint_{C^*} E(u') 
\frac{\theta^*\left(u-u'+c/2-P+1\right)}
{\theta^*(u-u'+c/2)}
\frac{\theta^*(1)}{\theta^*(P-1)}
\frac{dz'}{2\pi i z'},
\lb{Eplus}\\
&&F^+(u)
=a \oint_{C} F(u') 
\frac{\theta\left(u-u'+P+h-1\right)}{\theta(u-u')}
\frac{\theta (1)}{\theta(P+h-1)}
\frac{dz'}{2\pi i z'}. 
\lb{Fplus}
\ena
Here the contours are
\bea
C^* &:& |p^*q^c z|<|z'|<|q^cz|, 
\lb{C*}\\
C &:& |pz|<|z'|<|z|,
\lb{C}
\ena
and the constants $a,a^*$ are chosen to satisfy 
\begin{eqnarray*}
{a^* a \theta^*(1)\kappa\over q-q^{-1}} =1.
\end{eqnarray*}
\end{dfn}

We have to be careful about the ordering of $P$ and $E(u)$, $F(u)$ in 
\eqref{Eplus}--\eqref{Fplus}, since they do not commute. 
In fact we have the following commutation relations. 
\be
&&[K(u),P]=K(u),\quad [E(u),P]=2E(u),\quad [F(u),P]=0,
\\
&&[K(u),P+h]=K(u),\quad [E(u),P+h]=0,\quad [F(u),P+h]=2F(u).
\en
The specification of the contour 
\eqref{C} should be understood as an abbreviation of the prescription
``$C$ is a simple closed curve encircling the poles 
$z'=p^nz$ ($n\ge 1$) of the integrand, 
but not containing $z'=p^nz$ ($n\le 0$) inside''. 
Similarly for \eqref{C*}. 

The half currents \eqref{Eplus}--\eqref{Fplus} can also be written 
in terms of the Fourier modes of the elliptic currents 
\eqref{dress1}--\eqref{dress2},
\be
e(z,p)=\sum_{n\in\Z}e_nz^{-n},
\qquad 
f(z,p)=\sum_{n\in\Z}f_nz^{-n}.
\en
Substituting the Laurent expansion 
\be
\frac{\theta(u+s)}{\theta(u)\theta(s)}
=-\sum_{n\in\Z}\frac{1}{1-q^{-2s}p^n}z^{-n+\frac{s}{r}}
\qquad  (z=q^{2u})
\en
valid in the domain $1<|z|<|p^{-1}|$, we obtain 
\bea
&&E^+(u)=
e^{2Q}a^*\theta^*(1) 
\sum_{n\in\Z} e_n\frac{1}{1-q^{2(P-1)} p^{*n}}
(q^cz)^{-n-\Delta^*_{-P-1}+\Delta^*_{-P+1}},
\lb{Fou1}\\
&&
F^+(u)=
-a \theta(1)
 \sum_{n\in\Z} f_n\frac{1}{1-q^{-2(P+h-1)}p^{n}}
z^{-n+\Delta_{-P-h-1}-\Delta_{-P-h+1}}.
\lb{Fou2}
\ena
Readers who prefer the formal series language  may take 
\eqref{Fou1}, \eqref{Fou2} as the definition of the half currents. 

We remark that a change of contours leads to different definitions. 
For instance, we can define another pair of currents 
$E^-(u)$, $F^-(u)$ by the formulas \eqref{Eplus}-\eqref{Fplus}, 
with $C^*$, $C$ changed respectively to 
$C^*_-~:~|q^c z|<|z'|<|p^{*-1}q^cz|$ and 
$C_-~ : ~|z|<|z'|<|p^{-1}z|$. 
Then we have
\be
-a^*\theta^*(1) E(u)&\!\!=\!\!& E^+(u)-E^-(u),\\
-a\theta(1) F(u)&\!\!=\!\!&F^+(u)-F^-(u).
\en
This looks similar to the decomposition of the total currents to 
`positive' and `negative' parts in \cite{EF}. 
Notice however that in our case 
all the Fourier components $e_n,f_n$ 
appear in  $E^+(u)$, $F^+(u)$, and hence 
the `half' currents already generate the full algebra. 
For this reason we will not consider  $E^-(u)$, $F^-(u)$ and 
the analog of the $L^-$-operator in \cite{EF}. 

{}From the commutation relations \eqref{Ucom1}--\eqref{Ucom10} for the total 
currents, we can obtain the relations for the half currents. 
Recall the function $\rho(u)$ in \eqref{rho1}-\eqref{rho11}, 
which satisfies
\begin{eqnarray*}
&&\rho(0)=1,
\quad 
\rho(1)=\frac{\theta^*(1)}{\theta(1)},\\
&&
\rho(u)\rho(-u)=1,
\quad
\rho(u)\rho(u+1)=\frac{\theta^*(u+1)}{\theta^*(u)}
\frac{\theta(u)}{\theta(u+1)}.
\end{eqnarray*}

\begin{prop}\lb{prop:4.2}
Set $u=u_1-u_2$. 
Then the following commutation relations hold: 
\bea
\!\!\!\!\!\!\!\!\!\!&&K^+(u_1)K^+(u_2)=\rho(u)K^+(u_2)K^+(u_1),
\lb{hf1}\\
\!\!\!\!\!\!\!\!\!\!&&
K^+(u_1)E^+(u_2)K^+(u_1)^{-1}=E^+(u_2)\frac{\theta^*(1+u)}{\theta^*(u)}
-E^+(u_1){\theta^*\left(1\right)\over\theta^*\left(P\right)}
\frac{\theta^*\left(P+u\right)}{\theta^*(u)},
\lb{hf2}\\
\!\!\!\!\!\!\!\!\!\!&&
K^+(u_1)^{-1}F^+(u_2)K^+(u_1)
=
\frac{\theta(1+u)}{\theta(u)}F^+(u_2)
-
{\theta\left(1\right)\over
\theta\left(P+h\right)}
\frac{\theta\left(P+h-u\right)}{\theta(u)}F^+(u_1),
\lb{hf3}\\
\!\!\!\!\!\!\!\!\!\!&&
\frac{\theta^*(1-u)}{\theta^*(u)}E^+(u_1)E^+(u_2)
+
\frac{\theta^*(1+u)}{\theta^*(u)}E^+(u_2)E^+(u_1)
\lb{hf4}\\
\!\!\!\!\!\!\!\!\!\!&&
\qquad=
E^+(u_1)^2
{\theta^*\left(1\right) \over
\theta^*\left(P-2\right)}
{\theta^*\left(P-2+u\right) \over
\theta^*\left(u\right)}
+E^+(u_2)^2{\theta^*\left(1\right) \over
\theta^*\left(P-2\right)}
{\theta^*\left(P-2-u\right) \over
\theta^*\left(u\right)},
\nonumber\\
\!\!\!\!\!\!\!\!\!\!&&
\frac{\theta(1+u)}{\theta(u)}F^+(u_1)F^+(u_2)
+\frac{\theta(1-u)}{\theta(u)}F^+(u_2)F^+(u_1)
\lb{hf5} \\
\!\!\!\!\!\!\!\!\!\!&&=
F^+(u_1)^2
{\theta\left(1\right) \over
\theta\left(P+h-2\right)}
{\theta\left(P+h-2-u\right) \over
\theta\left(u\right)}
+F^+(u_2)^2
{\theta\left(1\right) \over
\theta\left(P+h-2\right)}
{\theta\left(P+h-2+u\right) \over
\theta\left(u\right)},
\nonumber\\
\!\!\!\!\!\!\!\!\!\!&&[E^+(u_1),F^+(u_2)]=
K^+(u_2-1)K^+(u_2)
\frac{\theta^*\left(P-1-u\right)}{\theta^*(u)}
\frac{\theta^*(1)}{\theta^*(P-1)} 
\nonumber\\
\!\!\!\!\!\!\!\!\!\!&&\qquad\qquad\qquad\quad\quad -
K^+(u_1)K^+(u_1-1)
\frac{\theta\left(P+h-1-u\right)}{\theta(u)}
\frac{\theta(1)}{\theta(P+h-1)}.
\lb{hf6}
\ena
\end{prop}

\proof 
These relations can be proven by reducing 
them to identities of theta functions. 
Let us show \eqref{hf2}. From the definition of the half currents  
\eqref{Eplus} and the commutation relation \eqref{Ucom2}, we have
\bea
&&K^+(u_1)E^+(u_2)K^+(u_1)^{-1}
\nonumber\\
&&=a^*\oint_{C^*}E(u')
\frac{\theta^*(u_1-u'+\frac{1-\rs}{2})}
{\theta^*(u_1-u'-\frac{1+\rs}{2})}
\frac{\theta^*\left(u_2-u'-P+1\right)}
{\theta^*(u_2-u')}
\frac{\theta^*(1)}{\theta^*(P-1)}
\frac{dz'}{2\pi i z'}.
\lb{hf2'}
\ena
Set 
\[
\eta_{s,t}(u)=
\frac{\theta(u+s)\theta(t)}{\theta(u)\theta(s)}.
\]
Then the following identity holds: 
\bea
&&
\frac{\theta(u_1+t)}{\theta(u_1)}
\eta_{s,t}(u_2)
=
\frac{\theta(u_1-u_2+t)}{\theta(u_1-u_2)}
\eta_{s+t,t}(u_2)
+
\eta_{s,t}(u_2-u_1)
\eta_{s+t,t}(u_1). 
\label{th1}
\ena
We obtain \eqref{hf2} by 
applying \eqref{th1} to the integrand of \eqref{hf2'} with 
the replacement 
$\theta(u) \rightarrow \theta^*(u)$, $a\rightarrow 1$,
$s\rightarrow-P+1$, $u\rightarrow u'$ and 
$u_1\rightarrow u_1-(1+\rs)/2$. 

Likewise, \eqref{hf4} leads to an equality between two-fold integrals. 
It can be shown by symmetrizing
the integration variables and applying the identity 
\be
&&
\psi(u_1,u_2;u_1',u_2')+
\psi(u_1,u_2;u_2',u_1')\times\frac{\theta(u_2'-u_1'+t)}{\theta(u_2'-u_1'-t)}
=( u_1 \longleftrightarrow u_2), 
\end{eqnarray*}
where
\bea
&&
\psi(u_1,u_2;u_1',u_2')=
\frac{\theta(u_1-u_2-t)}{\theta(u_1-u_2)}
\eta_{s+t,t}(u_1-u_1')
\eta_{s-t,t}(u_2-u_2')
\nonumber\\
&&
\qquad\qquad \qquad\qquad -
\eta_{s,t}(u_2-u_1)
\eta_{s+t,t}(u_1-u_1')\eta_{s-t,t}(u_1-u_2').
\label{th2}
\ena
The proofs of \eqref{hf3}, \eqref{hf5} are similar. 

Finally let us show \eqref{hf6}. 
Integrating the delta function in \eqref{Ucom6}, we obtain 
\be
&&(a^*a)^{-1}(q-q^{-1})[E^+(u_1),F^+(u_2)]
\\
&&\qquad =
\oint_{C_1} H^+(u'+\frac{c}{4})
\frac{\theta^*(u_1-u'-P+1)\theta^*(1)}{\theta^*(u_1-u')\theta^*(P-1)}
\frac{\theta(u_2-u'+P+h-1)\theta(1)}{\theta(u_2-u')\theta(P+h-1)}
\frac{dz'}{2\pi i z'}
\\
&&\qquad -
\oint_{C_2} H^-(u'-\frac{c}{4})
\frac{\theta^*(u_1-u'+c-P+1)\theta^*(1)}{\theta^*(u_1-u'+c)\theta^*(P-1)}
\frac{\theta(u_2-u'+P+h-1)\theta(1)}{\theta(u_2-u')\theta(P+h-1)}
\frac{dz'}{2\pi i z'}.
\en
Here the contours are 
\be
C_1 &:& |p^* z_1|,|pz_2|<|z'|<|z_1|,|z_2|, 
\\
C_2  &:& |p z_1|,|pz_2|<|z'|<|q^{2c}z_1|,|z_2|.
\en
Change variables $z'\rightarrow pz'$ in the second term 
and use the periodicity of $\theta(u)$ along with the 
relation $H^-(u'-c/4)=H^+(u'-r+c/4)$. We see that 
the integrand becomes the same as the first, whereas the contour becomes 
\be
C_2'  &:& |z_1|,|z_2|<|z'|<|p^{-1}q^{2c}z_1|,|p^{-1}z_2|.
\en
Picking the residues at $z'=z_1,z_2$ we find \eqref{hf6}. 
\qed

\subsection{Gau{\ss} decomposition}\lb{subsec:4.2}

Our next task is to 
rewrite the commutation relations \eqref{hf1}--\eqref{hf6}
into an `$RLL$'-form. 
Following the idea of the Gau{\ss} decomposition of Ding-Frenkel \cite{DF}, 
let us introduce the $L$-operator as follows. 

\begin{dfn}[$L$-operator]
We define the operator $\hL^+(u)\in \End(V)\otimes\Uqp\bigl(\slth\bigr)$ 
with $V=\C^2$, 
by 
\bea
&&\hL^+(u)=
\left(
\begin{array}{cc}
1 &F^+(u) \\
0 &1            \\
\end{array}
\right)
\left(
\begin{array}{cc}
K^+_1(u) & 0 \\
0             &K^+_2(u)\\
\end{array}
\right)
\left(
\begin{array}{cc}
1            &0      \\
E^+(u) &1      \\
\end{array}
\right),
\lb{Gauss}
\ena
where
\be
&&K^+_1(u)=K^+(u-1),\qquad 
K^+_2(u)=K^+(u)^{-1}.
\en
\end{dfn}
Note that 
\be
&&[\hL^+(u),h^{(1)}+h]=0,
\qquad
P \hL^+(u)=\hL^+(u)(P-h^{(1)}),
\en
where $h^{(1)}$ and $h$ mean 
$h\otimes 1$ and $1\otimes h\in \End(V)\otimes 
\Uqp\bigl(\slth\bigr)$, respectively. 

We also need the formula for the $R$ matrix
\eqref{tR} for $l=m=1$. 
With a further transformation of the form
\eqref{Rpr},  $R^+(u,s)=\tR^+_{11}(u,s)$ takes the form 
\bea
&&R^+(u,s)=\rho^+(u)
\left(
\begin{array}{cccc}
1 &                  &             & \\
  &b(u,s)      &c(u,s) & \\
  &\bar{c}(u,s)&\bar{b}(u,s) & \\
  &                  &             &1 \\
\end{array}
\right).
\lb{Rmat}
\ena
Here $\rho^+(u)$\footnote{This scalar factor 
differs from (3.16) in \cite{JKOS1},
see the remark at the end of section \ref{sec:2.1}.}
 is given in \eqref{rho11}, and 
\bea
&&b(u,s)=
\frac{\theta(s+1) \theta(s-1)  }{\theta(s)^2}
\frac{\theta(u)}{\theta(1+u)},
\quad 
c(u,s)=
\frac{\theta(1)}{\theta(s)}
\frac{\theta(s+u)}{\theta(1+u)},
\lb{Rmat1}\\
&&\bar{c}(u,s)=\frac{\theta(1)}{\theta(s)}
\frac{\theta(s-u)}{\theta(1+u)},
\quad \qquad\qquad\quad\; 
\bar{b}(u,s)=
\frac{\theta(u)}{\theta(1+u)}.
\lb{Rmat4}
\ena
Up to a scalar factor, 
this is the same $R$ matrix as eq.(93) in \cite{EF}. 

\begin{prop}\lb{prop:RLL1}
The relations \eqref{hf1}-\eqref{hf6} are 
equivalent to the following $RLL$ relation 
\bea
&&R^{+(12)}(u_1-u_2,P+h)
\hL^{+(1)}(u_1)\hL^{+(2)}(u_2)
\lb{RLL3}\\
&&\qquad =
\hL^{+(2)}(u_2)\hL^{+(1)}(u_1)
R^{*+(12)}(u_1-u_2,P).
\nn
\ena
\end{prop}
Proposition \ref{prop:RLL1} 
can be shown by a direct computation. 
(A little care has to be taken since 
the entries of $\hL^+(u)$ do not commute with those of $R^+(u,P)$.)
Since the calculation is tedious but straightforward, 
we omit the details. 

The above $RLL$-relation is equivalent to the dynamical 
$RLL$ relation \eqref{RLL2}. 
To see this, let us `strip off' the operator $e^Q$ from the half currents 
and define 
\begin{eqnarray*}
&&k^+_1(u,P)=
k(q^{r-1}z) \times
(q^{r-1}z)^{-{r-r^* \over 4 rr^*}(2P+1) +{1 \over 2 r} h}
=K^+(u-1)e^{-Q},\\
&&k^+_2(u,P)=
k(q^{r+1}z)^{-1} \times
(q^{r+1}z)^{{r-r^* \over4 rr^*}(2P-1) -{1 \over 2 r} h}
=K^+(u)^{-1}e^Q,
\\
&&e^+(u,P)=
 a^*\theta^*(1) 
\sum_{n\in\Z} e_n\frac{1}{1-q^{2P} p^{*n}}(q^cz)^{-n-\frac{P}{r^*}}
=e^{-Q}E^+(u)e^{-Q},
\\
&&
f^+(u,P)=
-a \theta(1)
 \sum_{n\in\Z} f_n\frac{1}{1-q^{-2(P+h-1)} p^{n}}z^{-n+\frac{P+h-1}{r}}
=F^+(u).
\end{eqnarray*}
These currents all commute with $P$. 
We can regard them as currents in $\uq$ having $P$ as a parameter 
($P$ plays the same role as $\la$ used in \cite{EF}). 
We set
\bea
L^+(u,P)
&=&\hL^+(u)
\left(\begin{array}{cc}
e^{-Q} &0\\
0 &e^{Q}\\
\end{array}\right)
\lb{Lplus}\\
&=&
\left(\begin{array}{cc}
1 &f^+(u,P) \\
0 &1            \\
\end{array}\right)
\left(\begin{array}{cc}
k_1^+(u,P) & 0 \\
0             &k_2^+(u,P)\\
\end{array}\right)
\left(\begin{array}{cc}
1            &0      \\
e^+(u,P) &1      \\
\end{array}\right).
\nn
\ena
Then Proposition \ref{prop:RLL1} is equivalently rephrased as follows. 
\begin{prop}\lb{prop:RLL2} 
The $L$-operator \eqref{Lplus} satisfies the dynamical $RLL$ relation 
\bea
&&R^{+(12)}(u_1-u_2,P+h)
L^{+(1)}(u_1,P)L^{+(2)}(u_2,P+h^{(1)})
\lb{RLL4}\\
&&\qquad =
L^{+(2)}(u_2,P)L^{+(1)}(u_1,P+h^{(2)})
R^{*+(12)}(u_1-u_2,P).
\nn
\ena
\end{prop}


\setcounter{section}{4}
\setcounter{equation}{0}
\section{Vertex operators}\lb{sec:5}

In this section we  shall study the VO's for $\Uqp=\Uqp\bigl(\slth\bigr)$ 
and compare them with those of Lukyanov and Pugai \cite{LukPug2}.

\subsection{Intertwining relations}\lb{subsec:5.1}

As in section \ref{sec:2}, we start with VO's 
in the sense of \eqref{Ph},\eqref{Ps} 
(or their modification by fractional powers \eqref{tPhi},\eqref{tPsi})
acting on some highest weight modules $\F_J$ over $\uq$, 
where $J$ is a label for the highest weight. 
Our main concern will be the Fock modules described in Appendix \ref{app:4}. 
However, for the general considerations till the end of the next subsection, 
we do not need the details of $\F_J$. 
For this purpose it is convenient to consider the VO's as acting on the sum 
of the Fock spaces $\F=\oplus_J \F_J$. 

We define the VO's for $\Uqp$ acting on the total Fock space 
$\hat{\F}=\oplus_\mu \F\otimes  e^{\mu Q}$ by
\bea
&&\hPhi_l(v)=\tPhi_l(v,P)~:~\hF\longrightarrow \hF\otimes V_{l,v},
\lb{Type1}\\
&&\hPsi^*_l(v)=\tPsi^*_l(v,P)e^{h\otimes Q }~:~
V_{l,v}\otimes \hF \longrightarrow \hF. 
\lb{Type2}
\ena
Here, just as in $\hLp(u)$, the $P$ in $\hPhi_l(v)$ and $\hPsi^*_l(v)$ 
is regarded as an operator on $\hF$. 
For the notation about the spin $l/2$ representation 
$V_{l,v}$, see Appendix \ref{app:3.1}.

The basic relations we are going to investigate are 
the dynamical intertwining relations in \eqref{int3} and \eqref{int4}.
We shall solve them by replacing $\tR^{+}_{1l}(u,P)$ and $\tL^{+}(u,P)$
with $R^{+}_{1l}(u,P)$ in \eqref{R1l1} and 
$L^{+}(u,P)$ in \eqref{Lplus} respectively.
Substituting \eqref{Lplus}, \eqref{Type1} and  \eqref{Type2}
into \eqref{int3} and \eqref{int4} and writing $u_1=u,\ u_2=v$,
we get the following for $\hPhi_l(v)$ and $\hPsi^*_l(v)$ :
\bea
&&\hPhi_l(v)\widehat{L}^+(u)=R^{+(13)}_{1l}(u-v,P+h)
\widehat{L}^+(u)\hPhi_l(v),
\lb{uqpint1}\\
&&
\widehat{L}^+(u)\widehat{\Psi}_l^*(v)
=\widehat{\Psi}_l^*(v)\widehat{L}^+(u)
R^{*+(12)}_{1l}(u-v,P-h^{(1)}-h^{(2)}).
\lb{uqpint2}
\ena

It should be noted that a natural coproduct for $\Uqp$ is not known, 
and hence the meaning of intertwining relations for it is not clear. 
Eqs. \eqref{uqpint1},\eqref{uqpint2} should be regarded 
as a compact way of writing the family of intertwining relations for 
$\Bqla(\slth)$. 
With this understanding, we shall sometimes refer to 
\eqref{uqpint1},\eqref{uqpint2} (somewhat loosely) 
as `intertwining relations for $\Uqp$'. 

Now using the explicit form of the $R$-matrix \eqref{R1l1}, 
let us write down the `intertwining relations' 
\eqref{uqpint1} and \eqref{uqpint2}. 
We shall write the entries of $\widehat{L}^+(u)$ as 
\be
&&\hLp(u)=\left(\matrix{\hLp_{++}(u)&\hLp_{+-}(u)\cr
\hLp_{-+}(u)&\hLp_{--}(u)\cr
}\right).
\en
According to the Gau{\ss} decomposition \eqref{Gauss}, we have
\be
&&\widehat{L}^+_{++}(u)=K^+(u-1)+F^+(u)K^+(u)^{-1}E^+(u),\quad
\widehat{L}^+_{+-}(u)=F^+(u)K^+(u)^{-1},
\\
&&\widehat{L}^+_{-+}(u)=K^+(u)^{-1}E^+(u),\qquad\qquad\qquad\qquad\quad\;\;
\widehat{L}^+_{--}(u)=K^+(u)^{-1}.
\en
Define the components of VO's by
\be
&&\hPhi_l(v-\frac{1}{2})=\sum_{m=0}^{l}\Phi_{l,m}(v)\otimes v_m^l,
\\
&&\hPsis_l(v-\frac{c+1}{2})\left(v_m^l\otimes \cdot \right)=\Psi^*_{l,m}(v).
\en
For brevity we set 
$\tilde{\varphi}_l(u)=\varphi_l(u+\frac{1}{2}),\ 
{\tilde{\varphi}^*}_l(u)=\varphi^*_l(u+\frac{c+1}{2})$ and  
$\La=P+h$. 
Here $\varphi_l(u)$ is given in \eqref{varph} and 
$\varphi^*_l(u)=\varphi_l(u)\bigl|_{r\to r^*}$.
In these notations, \eqref{uqpint1} reads as follows: 
\bea 
&&\tilde{\varphi}_l(u-v)\Phi_{l,m}(v)\hLp_{+\pm}(u)
=\nonumber\\
&&\qquad\qquad
-\frac{\theta(u-v+\frac{l}{2}-m+1)\theta(\La+l-m+1)\theta(\La-m)}
{\theta(\La)\theta(\La+l-2m+1)} \hLp_{+\pm}(u)\Phi_{l,m}(v)
\nonumber\\
&&\qquad \qquad-
\frac{\theta(u-v+\La+\frac{l}{2}-m+1)\theta(m)}{\theta(\La+l-2m+1)}
\hLp_{-\pm}(u)\Phi_{l,m-1}(v),
\label{ppm1}\\
&&\tilde{\varphi}_l(u-v)\Phi_{l,m}(v)\hLp_{-\pm}(u)
=-\theta(u-v-\tfrac{l}{2}+m+1)\hLp_{-\pm}(u)\Phi_{l,m}(v)
\nonumber\\
&&\qquad\qquad 
+\frac{\theta(u-v-\La-\frac{l}{2}+m+1)\theta(l-m)}{\theta(\La)}
\hLp_{+\pm}(u)\Phi_{l,m+1}(v)\label{mpm1}. 
\ena
Similarly \eqref{uqpint2} takes the form 
\bea 
&&{\tilde{\varphi}^*}_l(u-v)\hLp_{\pm +}(u)\Psi^*_{l,m}(v)
=\nonumber\\
&&\qquad\qquad
-\frac{\theta^*(u-v+\frac{l+c}{2}-m+1)\theta^*(P-l+m-1)\theta^*(P+m)}
{\theta^*(P)\theta^*(P-l+2m-1)}\Psi^*_{l,m}(v) \hLp_{\pm +}(u)
\nonumber\\
&&\qquad \qquad+
\frac{\theta^*(u-v-P+\frac{l+c}{2}-m+1)\theta^*(l-m+1)}
{\theta^*(P-l+2m-1)}\Psi^*_{l,m-1}(v)\hLp_{\pm -}(u),
\label{pmp2}\\
&&{\tilde{\varphi}^*}_l(u-v)\hLp_{\pm -}(u)\Psi^*_{l,m}(v)
=
-\theta^*(u-v-\tfrac{l-c}{2}+m+1)\Psi^*_{l,m}(v)\hLp_{\pm -}(u)
\nonumber\\
&&\qquad\qquad 
-\frac{\theta^*(u-v+P-\frac{l-c}{2}+m+1)\theta^*(m+1)}{\theta^*(P)}
\Psi^*_{l,m+1}(v)\hLp_{\pm +}(u).\label{pmm2} 
\ena

Let us investigate the relations \eqref{ppm1} and \eqref{mpm1} in detail. 
For the highest component $\Phi_{l,l}(v)$, we can immediately obtain 
{}from \eqref{mpm1} the relations
\bea
&&
\tilde{\varphi}_l(u-v)\Phi_{l,l}(v)K^+(u)^{-1}=-
\theta(u-v+\frac{l}{2}+1)K^+(u)^{-1}\Phi_{l,l}(v),
\label{intpk}\\
&&
\Phi_{l,l}(v)E^+(u)=E^+(u)\Phi_{l,l}(v)\label{intpe}.
\ena
Notice that $u=v-\frac{l}{2}-1$ is a zero of $\tilde{\varphi}_l(u-v)$. 
Suppose that in the relation \eqref{ppm1}
the product of $\Phi_{l,m}(v)$ and $F^+(u)$ has no pole at this point. 
Then we obtain a relation which determines the components 
of VO's recursively:
\bea
\Phi_{l,m-1}(v)=F^+(v-\frac{l}{2})\frac{\theta(\La+l-m)}{\theta(\La)}
\Phi_{l,m}(v)\qquad (m=0,1,..,l)\label{lowerp}.
\ena
We will show below (Proposition \ref{prop:5.1}) that this assumption is 
satisfied in the free field realization of $\Uqp$. 
Substituting \eqref{intpk}, \eqref{lowerp} into \eqref{ppm1} and using 
Riemann's theta identity, 
we find the following relation as the sufficient condition  
for \eqref{ppm1} with the choice of the lower sign:
\bea
&&
\Phi_{l,l}(v)F(u)=\frac{\theta(u-v-\frac{l}{2})}{\theta(u-v+\frac{l}{2})}
F(u)\Phi_{l,l}(v).
\label{intpf}
\ena

By \eqref{intpk}-\eqref{intpf}, 
the remaining relations in \eqref{ppm1}-\eqref{mpm1} are
reduced to those involving the highest component $\Phi_{l,l}(v)$ and  
$K^+(u),E^+(u),F^+(u)$.  In order to ensure the  
existence of VO, we need to verify them.
For level one ($c=1$), we have verified that
they are consequences of Proposition \ref{prop:4.2}.  
In general, such a direct check seems complicated, and it would be  
better to invoke the fusion procedure.
We do not go into this issue any further. Note however that, had we  
known the equivalence with the quasi-Hopf construction, there would  
be no need for the check because the existence of VO is clear in the latter  
context.

Similarly, the intertwining relations \eqref{pmp2} and \eqref{pmm2} 
for the type II vertex operator lead to the following relations as the 
sufficient condition for the highest component:
\bea
&&{\tilde{\varphi}^*}_l(u-v)K^+(u)^{-1}\Psi^*_{l,l}(v)
=-\Psi^*_{l,l}(v)K^+(u)^{-1}\theta^*(u-v+\frac{l+c}{2}+1),
\label{intpsik}\\
&&\Psi^*_{l,l}(v)E(u)=E(u)\Psi^*_{l,l}(v)
\frac{\theta^*(u-v+\frac{l}{2})}
{\theta^*(u-v-\frac{l}{2})},
\label{intpsie}\\
&&\Psi^*_{l,l}(v)F^+(u)=F^+(u)\Psi^*_{l,l}(v),\label{intpsif}
\ena
and the relation for the lower component
\begin{equation}
\Psi^*_{l,m-1}(v)=\Psi^*_{l,m}(v)E^+(v-\frac{l+c}{2}-r^*)
\frac{\theta^*(m)\theta^*(P-l+m-2)}
{\theta^*(l-m+1)\theta^*(P-2)}\quad (m=0,1,..,l).
\label{lowerpsi}
\end{equation}
We remark that in the derivation of \eqref{lowerpsi}, 
we took $u=v-\frac{l+c}{2}-1-r^*$ as 
a zero of $\tilde\varphi^*_l(u-v)$. If we chose
a zero without the shift by $r^*$, we would have an extra term 
in the RHS of \eqref{intpsie}.
The shift of $u$ by $r^*$ in $E^+(u)$ yields a change of the contour
in \eqref{Eplus}. 
For example, 
we have
\be
&&E^+(v-\frac{l+c}{2}-r^*)
=a^* \oint_{\tilde{C}^*} E(u') 
\frac{\theta^*\left(v-u'-l/2-P+1\right)}
{\theta^*(u-u'-l/2)}
\frac{\theta^*(1)}{\theta^*(P-1)}
\frac{dz'}{2\pi i z'}
\en
with the contour being
\bea
\tilde{C}^* &:& |q^{-l} z|<|z'|<|p^{*-1}q^{-l}z|.
\label{tildeCstar}
\ena

\subsection{`Twisted' intertwining relations}\lb{subsec:5.2}

In order to compare 
these results with those in \cite{LukPug2}, \cite{MW96}
and \cite{Konno}, we need a further modification by signs. 
In \cite{KLP1}, the modified VO's are called `twisted' intertwiners
\footnote{This terminology is not to be confused with the `twisting'
in the sense of quasi-Hopf algebras.}. 

In the case of $U_q(\slth)$, the twisted type I $\pPhi_V(v)$ and the type II 
$\pPsi_V(v)$ VO's are the 
intertwiners  of the same type as $\Phi_V(v)$ and
$\Psi^*_V(v)$,  but
satisfying the following intertwining relations twisted with signs:
\bea
&&\pPhi_V(v)\iota(a) = \Delta(a)\pPhi_V(v),\\
&&\pPsi_V(v)\Delta(a) = \iota(a) \pPsi_V(v),\qquad \forall a\in U_q(\slth), 
\ena
where $\iota$ denotes the involution of $ U_q(\slth)$,
\be
\iota(x^{\pm}_n)=-x^{\pm}_n,\qquad \iota(a_n)=a_n.
\en



Analogously to the procedures \eqref{Ph}, \eqref{Ps}, \eqref{Type1} and 
\eqref{Type2}, we can define the `twisted' VO's 
${\thPhi}_l(v)$  and ${\thPsis}_l(v)$ for $\Uqp$ as the operators of 
the same type as \eqref{Type1} and \eqref{Type2} but satisfying the following 
`twisted' intertwining relations.
\bea
&&{\thPhi}_l(v)\thL(u)=R^{+(13)}(u-v,P+h)
\hLp(u){\thPhi}_l(v),
\lb{inttilde1}\\
&&
\thL(u){\thPsis}_l(v)
={\thPsis}_l(v)\hLp(u)R^{*+(12)}(u-v,P-h^{(1)}-h^{(2)}).
\lb{inttilde2}
\ena
Here the `twisted' $L$-operator $\thL(u)$ has the following components: 
$$
\thL_{\pm\pm}(u)=\widehat{{L}}^+_{\pm\pm}(u),\qquad
\thL_{\pm\mp}(u)=-\widehat{{L}}^+_{\pm\mp}(u).
$$ 
Defining the components of the `twisted' VO's as
\be
&&\widehat{\Phi}'_l(v-\frac{1}{2})=\sum_{m=0}^{l}\pPhi_{l,m}(v)\otimes v_m^l,
\\
&&{\thPsis}_l(v-\frac{c+1}{2})\left(v_m^l \otimes \cdot \right)=\pPsi_{l,m}(v),
\en
we obtain the `twisted' counterpart of \eqref{intpk}-\eqref{lowerpsi} as 
follows:
\bea
&&
\tilde{\varphi}_l(u-v)\pPhi_{l,l}(v)K^+(u)^{-1}=
-\theta(u-v+\frac{l}{2}+1)K^+(u)^{-1}\pPhi_{l,l}(v),
\label{tintpk}\\
&&
\pPhi_{l,l}(v)E^+(u)=-E^+(u)\pPhi_{l,l}(v)\label{tintpe},\\
&&
\pPhi_{l,l}(v)F(u)=-\frac{\theta(u-v-\frac{l}{2})}{\theta(u-v+\frac{l}{2})}
F(u)\pPhi_{l,l}(v),\label{tintpf}\\
&&\pPhi_{l,m-1}(v)=F^+(v-\frac{l}{2})\frac{\theta(\La+l-m)}{\theta(\La)}
\pPhi_{l,m}(v)\label{tlowerp},
\ena
for type I, and 
\bea
&&{\tilde{\varphi}^*}_{l}(u-v)K^+(u)^{-1}\pPsi_{l,l}(v)=
-\pPsi_{l,l}(v)K^+(u)^{-1}
\theta^*(u-v+\frac{l+c}{2}+1),\label{tintpsik}\\
&&\pPsi_{l,l}(v)E(u)=-E(u)\pPsi_{l,l}(v)
\frac{\theta^*(u-v+\frac{l}{2})}
{\theta^*(u-v-\frac{l}{2})},\label{tintpsie}\\
&&\pPsi_{l,l}(v)F^+(u)=-F^+(u)\pPsi_{l,l}(v),\label{tintpsif}\\
&&\pPsi_{l,m-1}(v)=\pPsi_{l,m}(v)E^+(v-\frac{l+c}{2}-r^*)
\frac{\theta^*(m)\theta^*(P-l+m-2)}
{\theta^*(l-m+1)\theta^*(P-2)},\label{tlowerpsi}
\ena
for type II. (See the remark below \eqref{lowerpsi}.)

\subsection{Free field realization}\lb{subsec:5.3}

$\Uqp\bigl(\slth\bigr)$ admits a free field representation for arbitrary level 
$c=k (\not=0,-2)$ \cite{Konno}(see Appendix D).
Using this, we obtain a realization of the VO's. 
We consider below the `twisted' VO's. 
The non-twisted ones are obtained in a similar way. 

\begin{prop}\lb{prop:5.1}
Using the notation in Appendix D, 
the intertwining relation \eqref{inttilde1} has the following solution:
\bea
\!\!\!&&\Phi'_{l,l}(v)=\phi_{k-l,-(k-l)}(w):
\EXP{-\phi_0'(l;2,k|w)}:\qquad (l=0,1,..,k),\label{ttypeivo}\\
\!\!\!&&\Phi'_{l,m}(v)=\left(\prod_{j=1}^{l-m} a\oint_{C_j}
\frac{d z_j}{2\pi i z_j}
\right)\prod_{j=1}^{l-m}\frac{\theta(v-u_j-\frac{l}{2}+\La-1+2j)
\theta(1)\theta(\La+l-m-1+2j)}{\theta(v-u_j-\frac{l}{2})
\theta(\La-1+2j)\theta(\La+2j)}
\nonumber\\
\!\!\!&&\qquad\quad\quad\times
\left(\prod_{j=1}^{l-m}F(u_j)\right)\Phi'_{l,l}(v)
\qquad (m=0,1,..,l),\label{ttypeilower}
\ena
where  $w=q^{2v},\ z_j=q^{2u_j}$, $\Lambda=P+h$ and 
\bea
&&\phi_{l,\pm l}(w)=\,:\EXP{-\phi_2\Bigl(\pm l;2,k|w;\pm\frac{k}{2}\Bigr)-
\phi_1\Bigl(l;2,k+2|w;\pm\frac{k+2}{2}\Bigr)}:.
\ena 
Similarly \eqref{inttilde2} has the solution
\bea
&&{\Psi^*}'_{l,l}(v)=\phi_{k-l,k-l}(w):
\EXP{\phi_0(l;2,k|z)}:\qquad(l=0,1,2,..,k),\label{ttypeiivo}\\
&&{\Psi^*}'_{l,m}(v)=\left(\prod_{j=1}^{l-m} a^*\oint_{\tilde{C}_j^*}
\frac{d z_j}{2\pi i z_j}
\right){\Psi^*}'_{l,l}(v)\left(\prod_{j=1}^{l-m}E(u_j)\right)
\nonumber\\
&&\qquad\times
\prod_{j=1}^{l-m}\frac{\theta^*(v-u_j-\frac{l}{2}-P+1+2(l-m-j))
\theta^*(1)\theta^*(P-j-1)\theta^*(l-j+1)}
{\theta^*(v-u_j-\frac{l}{2})
\theta^*(P-1-2(l-m-j))\theta^*(P-2-2(l-m-j))\theta^*(j)}
\nonumber\\
&&\qquad \qquad\qquad\qquad\qquad\qquad\qquad \qquad 
(m=0,1,..,l).\label{ttypeiilower}
\ena
The contours $C_j$ and $\tilde{C}^*_j\ (j=1,2,..,l-m)$ are taken as
\bea
&&C_{l-m}\ :\ |q^lw|,|pq^{-l}w|\ <\ |z_{l-m}|\ <\ |q^{-l}w|,\label{cont1}\\
&&C_{j}\ :\ |q^lw|,|pq^{-l}w|, |q^{-2}z_{j+1}|\ <\ |z_j|\ <\ |q^{-l}w|
\quad (j=1,2,..,l-m-1),\label{cont2}\\
&&\tilde{C}^*_{l-m}\ :\ |q^{-l}w|\ <\ |z_{l-m}|\ <\ |p^{*-1}q^{-l}w|,
|q^{l}w|,\label{conts1}\\
&&\tilde{C}^*_{j}\ :\ |q^{-l}w|, |q^{2}z_{j+1}|\ <\ |z_j|\ <\ 
|p^{*-1}q^{-l}w|, |q^{l}w|
\ (j=1,2,..,l-m-1)\label{conts2}.
\ena
\end{prop}

\noindent {\it Sketch of proof.}\quad 
First of all, 
we should note that the relations \eqref{tintpk}-\eqref{tintpf} and 
\eqref{tintpsik}-\eqref{tintpsif} coincide with those in Proposition 4.4 
in \cite{Konno} if we identify
 $\Phi'_{l,l}(u)$ and ${\Psi^*}'_{l,l}(u)$
with $\tilde\Phi^{(l)}_l(u)$ and $\tilde\Psi^{(l)*}_l(u)$ in \cite{Konno}, 
respectively. 
Using the bosons $a_{0,m}, a_{1,m}, a_{2,m}$ defined 
in Appendix D, the expressions \eqref{ttypeivo} and
\eqref{ttypeiivo} are the unique solutions of them up to a scalar.

Next, the expressions \eqref{ttypeilower} and \eqref{ttypeiilower}
are the direct consequence of \eqref{tlowerp} and \eqref{tlowerpsi}.
The contours are determined from \eqref{C} with the 
replacement $z\to q^{-l}w$, 
\eqref{tildeCstar} and the 
following OPE's derived from \eqref{ttypeivo}, \eqref{ttypeiivo}
and Proposition D.3:
\bea
&&\Phi'_{l,l}(v)F(u)=w^{\frac{l}{r}-1}
\frac{(pq^{-l}z/w;p)_{\infty}}{(q^{l}z/w;p)_{\infty}}:\Phi'_{l,l}(v)F(u):,
\label{opephif}\\
&&F(u)\Phi'_{l,l}(v)=z^{\frac{l}{r}-1}
\frac{(pq^{-l}w/z;p)_{\infty}}{(q^{l}w/z;p)_{\infty}}:\Phi'_{l,l}(v)F(u):,
\label{opefphi}\\
&&\Psi^{*'}_{l,l}(v)E(u)=w^{-\frac{l}{r-k}-1}
\frac{(p^*q^{l}z/w;p^*)_{\infty}}{(q^{-l}z/w;p^*)_{\infty}}:E(u)
\Psi^{*'}_{l,l}(v): ,\label{opepsie}\\
&&E(u)\Psi^{*'}_{l,l}(v)=z^{-\frac{l}{r-k}-1}
\frac{(p^*q^{l}w/z;p^*)_{\infty}}{(q^{-l}w/z;p^*)_{\infty}}:E(u)
\Psi^{*'}_{l,l}(v): ,\label{opeepsi}\\
&&E(u)E(u')=z^{\frac{2}{r-k}}
\frac{(p^*q^{-2}z'/z;p^*)_{\infty}}{(q^{2}z'/z;p^*)_{\infty}}
:e^{-\phi_0(k|z)}e^{-\phi_0(k|z')}: \nonumber\\
&&\qquad\qquad\times\frac{1}{(q-q^{-1})^2}
\Bigl(q(1-\frac{z'}{z})(:\Psi^-_I(z)\Psi^-_I(z'):+ 
:\Psi^-_{II}(z)\Psi^-_{II}(z'):\nonumber\\
&&\qquad\qquad\qquad-(1-q^{-2}\frac{z'}{z})(q :\Psi^-_I(z)\Psi^-_{II}(z'):
+ q^{-1}:\Psi^-_{II}(z)\Psi^-_I(z'): )
\Bigr),\label{opeee}\\
&&F(u)F(u')=z^{-\frac{2}{r}}
\frac{(pq^{2}z'/z;p)_{\infty}}{(q^{-2}z'/z;p)_{\infty}}
:e^{\phi'_0(k|z)}e^{\phi'_0(k|z')}: \nonumber\\
&&\qquad\qquad\times\frac{1}{(q-q^{-1})^2}
\Bigl(q^{-1}(1-\frac{z'}{z})(:\Psi^+_I(z)\Psi^+_I(z'):
+ :\Psi^+_{II}(z)\Psi^+_{II}(z'):\nonumber\\
&&\qquad\qquad\qquad
-(1-q^{2}\frac{z'}{z})(q^{-1} :\Psi^+_I(z)\Psi^+_{II}(z'):
+ q:\Psi^+_{II}(z)\Psi^+_I(z'): )
\Bigr),\label{opeff}
\ena
where $\Psi^\pm_{I,II}(z)$ are given in \eqref{pfcurrents} and 
\eqref{fullpf}, and  $z=q^{2u},\ z'=q^{2u'}$, $w=q^{2v}$.

In the derivation of \eqref{tlowerp} and \eqref{tlowerpsi}, 
we made an assumption that
no counter poles to the zeros of $\tilde{\varphi}_l(u-v)$ and
${\tilde{\varphi}^*}_l(u-v)$ appear from the OPE's $\Phi'_{l,m}(v)F^+(u)$ 
and $E^+(u)\Psi^{*'}_{l,m}(v)$.
The verification of this assumption is not so hard.     
Substitute the OPE's \eqref{opephif}- \eqref{opeff}
into the products $\Phi'_{l,m}(v)F^+(u)$ 
and $E^+(u)\Psi^{*'}_{l,m}(v)$, we can show 
that such counter poles do not appear.

We have also checked at level one ($k=1$) that, 
upon elimination of ${\Phi'}_{1,m}(v)$ and $\Psi^{*'}_{1,m}(v)$, 
the rest of the intertwining relations are consequences of the 
commutation relations for the half currents in Proposition 4.2. 
As we mentioned already, such a direct check is  difficult for higher levels.  
However, modulo the assumption about 
the equivalence with the quasi-Hopf construction 
(where the existence of the VO's is known),  
the expressions \eqref{ttypeivo}-\eqref{ttypeiilower} are unique up to 
a scalar factor and a choice of an equivalent set of three bosons. 
\qed

\subsection{Level one case $(k=1)$}\lb{subsec:5.4}

In this case, writing $\Phi_-(v)=\Phi'_{1,1}(v),\ \Phi_+(v)=\Phi'_{1,0}(v),
\ \Psi^*_-(v)=\Psi^{*'}_{1,1}(v)$ and $\Psi^*_+(v)=\Psi^{*'}_{1,0}(v)$, 
we have from Proposition \ref{prop:5.1}, 
\bea
&&\Phi_-(v)=\,:\EXP{-\phi_0'(2|w)}:,\label{level1pm}\\
&&\Phi_{+}(v)= a\oint_{C}\frac{d z'}{2\pi i z'}
z'^{-\frac{r-1}{r}}\frac{(pq^{-1}w/z';p)_{\infty}}{(q w/z';p)_{\infty}}
\frac{\theta(v-u'-\frac{1}{2}+\La+1)
\theta(1)}{\theta(v-u'-\frac{1}{2})
\theta(\La+1)}
:F(u')\Phi_-(v):,\nonumber\\
\label{lpip}\\
&&{\Psi^*}_-(v)=\,:\EXP{\phi_0(2|z)}:,\label{lpiim}\\
&&{\Psi^*}_+(v)=a^*\oint_{\tilde{C}^*}\frac{d z'}{2\pi i z'}
:\Psi^*_-(v)E(u'): w^{-\frac{r}{r-1}}
\frac{(p^*q z'/w;p^*)_{\infty}}{(q^{-1} z'/w;p^*)_{\infty}}
\frac{\theta^*(v-u'-\frac{1}{2}-P+1)
\theta^*(1)}{\theta^*(v-u'-\frac{1}{2})\theta^*(P-1)},\nonumber\\
\label{level1psip}
\ena
where the contours $C$ and $\tilde{C}^*$ are given by \eqref{cont1} and 
\eqref{conts1} letting $l=1,\ m=0$.
Since the level one parafermion theory is a trivial theory, one can neglect 
the parafermion currents in $E(u)$ and $F(u)$. 
Then, the expressions \eqref{level1pm}-\eqref{level1psip} 
agree with the results in \cite{LukPug2} and \cite{MW96}.

\Remark Our notation here is related to those of \cite{MW96} as follows:
\be
&&x=-q, \quad \alpha_m=\frac{1}{[2m]}a_{0,m}, \quad
\beta_n=\frac{1}{[2m]}a'_{0,m},
\\
&&Q=Q_0, \quad P=P_0,
\\
&&L-1=P\ ({\rm in\ \eqref{level1psip}})=\hat{\Pi}+1, 
\quad K-1=\La=\hat{\Pi}'+1,
\\ 
&&A(z)=F(u), \quad B(z)=E(u).
\en


\setcounter{section}{5}
\setcounter{equation}{0}
\section{Discussions} \lb{sec:6}

\subsection{Classical limit}\lb{sec:6.1}

We have not been able so far to identify 
the $L$-operator constructed in section \ref{sec:4} 
with the one obtained by a quasi-Hopf twist \cite{JKOS1}. 
In this subsection, we study the classical limit 
of the elliptic algebra ${\cal B}_{q,\la}\bigl(\slth\bigr)$ 
in the $RLL$ formulation, 
and compare it with the half currents \eqref{Fou1}, \eqref{Fou2} 
in the quantum case. 

Let $\goth{a}=\slt$, with the standard generators $e,f,h$. 
Let $(~,~)$ be the invariant inner product normalized as 
 $(h,h)=2$, $(e,f)=(f,e)=1$. 
We consider the homogeneous realization of 
the affine Lie algebra $\goth{g}=\widehat{\slt}$, 
\be
&&\goth{g}=\mbox{span}\{ e_n,f_n,h_n \  (n\in \Z), d,c\}, 
\\
&&[x_m,y_n]=[x,y]_{m+n}+m(x,y)\delta_{m+n,0}c,\\
&&[d,x_m]=mx_m, \\
&&c:\hbox{ central }. 
\en
In what follows, we identify $h_0\in\goth{g}$ with $h\in\goth{a}$. 

Let us recall the notion of a quasi-Lie bialgebra \cite{QHA}
which is the classical counterpart of a quasi-Hopf algebra.
By definition, it is a triple $(\goth{g},\delta,\varphi)$ 
consisting of a Lie algebra $\goth{g}$, a $1$-cocycle (cobracket) 
$\delta:\goth{g}\rightarrow \wedge^2\goth{g}$ and a tensor 
$\varphi\in\wedge^3\goth{g}$, satisfying 
\bea
&&\frac{1}{2}\mbox{Alt}(\delta\otimes 1)\delta(x)
=[x^{(1)}+x^{(2)}+x^{(3)},\varphi],
\lb{qLalg1}\\
&&\mbox{Alt}(\delta\otimes 1\otimes 1)\varphi=0.
\lb{qLalg2}
\ena
Here the symbol Alt stands for skew-symmetrization. 
In the case of ${\cal B}_{q,\la}\bigl(\slth\bigr)$, the corresponding 
quasi-Lie bialgebra structure on $\goth{g}$ is 
given as follows \cite{Fron1}: 
\bea
&&\delta (x)= [x^{(1)}+x^{(2)},r],
\lb{cobra}\\
&&\varphi=-2\left(D^{(1)}r^{(23)}
-D^{(2)}r^{(13)}+D^{(3)}r^{(12)}\right).
\lb{SOSphi}
\ena
Here $r$ denotes the classical $r$ matrix 
\bea
&&r=\frac{1}{2}h\otimes h+
\sum_{n \neq 0}\frac{1}{1-p^n}h_n\otimes h_{-n}
+2\sum_{n\in\Z} \frac{1}{1-w p^n}e_n\otimes f_{-n}
\nn\\
&&\qquad 
+2\sum_{n\in\Z} \frac{1}{1-w^{-1}p^n}f_n\otimes e_{-n}
+c\otimes d+ d\otimes c,
\lb{SOSr}
\ena
with $p, w$ being parameters
(having the same meaning as in the body of the text),  
and we have set 
\be
D^{(i)}r^{(jk)}=
\left(c^{(i)}p\frac{\partial}{\partial p}
+h^{(i)}w\frac{\partial}{\partial w}\right)r^{(jk)}.
\en

Let $\rho_z:\goth{g}'\rightarrow \goth{a}\otimes\C[z,z^{-1}]$ 
($\goth{g}'=[\goth{g},\goth{g}]$) be the evaluation morphism given by 
$\rho_z(x_n)=xz^n$ ($x=e,f,h$), $\rho_z(c)=0$. 
Setting $r'=r-c\otimes d-d\otimes c$, we define 
\be
&&\cL^+(z)=\left(\rho_z\otimes\id\right)r',
\qquad
r(z)=\left(\id\otimes \rho_{1}\right)\cL^+(z).
\en
These are formal series with values in $\goth{a}\otimes\goth{g}$ and 
$\goth{a}\otimes\goth{a}$ , respectively. 
Explicitly we have 
\bea
\cL^+(z)&\!\!=\!\!&
h\otimes\left(\frac{1}{2}h+\sum_{n\neq 0}\frac{1}{1-p^n}z^n h_{-n}\right)
\nn\\
&&
+e\otimes\left(\sum_{n\in\Z}\frac{2}{1-w p^n}z^n f_{-n}\right)
+f\otimes\left(\sum_{n\in\Z}\frac{2}{1-w^{-1} p^n}z^n e_{-n}\right),
\lb{L+}\\
r(z)&\!\!=\!\!&
h\otimes h\left(\frac{1}{2}+\sum_{n\neq 0}\frac{1}{1-p^n}z^{n}\right)
\nn\\
&&
+e\otimes f \left(\sum_{n\in\Z}\frac{2}{1-w p^n}z^{n}\right)
+f\otimes e \left(\sum_{n\in\Z}\frac{2}{1-w^{-1} p^n}z^{n}\right).
\lb{clr}
\ena
Up to a change of `gauge' and the extra zero-mode operators $P,Q$, 
\eqref{L+} agrees with the classical limit of the $L$ operator based on 
the half currents \eqref{Fou1}, \eqref{Fou2}. 

According to Drinfeld \cite{QHA}, there is a bijective correspondence between 
quasi-Lie bialgebras and Manin pairs. 
This means the following. 
Let $D\goth{g}=\goth{g}\oplus\goth{g}^*$ be 
the direct sum with the dual vector space $\goth{g}^*$. 
Equip it with the inner product $(~,~)$ requiring that it 
vanishes on $\goth{g}\times\goth{g}$, $\goth{g}^*\times\goth{g}^*$ and 
coincides with the canonical paring on $\goth{g}\times\goth{g}^*$. 
Then, for  $(\goth{g},\delta,\varphi)$ a quasi-Lie bialgebra, 
 $D\goth{g}$ is endowed with a unique Lie algebra structure 
(the classical double) such that 
$(~,~)$ is an invariant inner product for $D\goth{g}$.
Moreover 
the correspondence $(\goth{g},\delta,\varphi)\leftrightarrow D\goth{g}$ 
is one-to-one. 

In the present case, let us take the dual basis 
$e_n^*,f_n^*,h_n^*$  ($n\in \Z$), $d^*,c^*$ of  
$\goth{g}^*$, with the dual pairing given by 
\be
&&\br{x_m,y_n^*}=\br{x,y}\delta_{m+n,0}, \quad
\br{d,c^*}=1, \quad \br{c,d^*}=1, \quad 
\mbox{ others }=0.
\en
Set\footnote{The $\cL^-(z)$ is a generating series in $\goth{g}^*$
and is independent of $\cL^+(z)$. 
It should not be confused with the classical limit of the $L^-(z)$ operator 
in \cite{JKOS1} which has a simple relation with $L^+(z)$.}
\bea 
&&\cL^-(z)=\frac{1}{2}h\otimes h^-(z)+e\otimes f^-(z)+f\otimes e^-(z),
\lb{Lm}
\\
&& x^-(z)=\sum_{n\in\Z} x^*_n z^{-n}
\qquad (x=e,f,h).
\nn
\ena
Then the dual pairing takes the form 
\be
\br{\cL^{+(1)}(z_1),\cL^{-(2)}(z_2)}=r^{(12)}(z_1/z_2).
\en

With the above notation, 
the Lie algebra structure of $D\goth{g}$ can be described as follows. 
\be
&&[\cL^{\pm(1)}(z_1),\cL^{\pm(2)}(z_2)]
=-[r^{(12)}(z_1/z_2),\cL^{\pm(1)}(z_1)+\cL^{\pm(2)}(z_2)]
\\
&&\quad \pm 2w\frac{\partial }{\partial w}
\left(h^{(1)}\cL^{+(2)}(z_2)-h^{(2)}\cL^{+(1)}(z_1)
+h r^{(12)}(z_1/z_2)\right)
\pm 2c \,p\frac{\partial }{\partial p}r^{(12)}(z_1/z_2),
\\
&&[\cL^{+(1)}(z_1),\cL^{-(2)}(z_2)]=
-[r^{(12)}(z_1/z_2),\cL^{+(1)}(z_1)+\cL^{-(2)}(z_2)]
\\
&&\quad +2 w\frac{\partial }{\partial w}
\left(h^{(1)}\cL^{+(2)}(z_2)-h^{(2)}\cL^{+(1)}(z_1)
+h r^{(12)}(z_1/z_2)\right)
\\
&&\quad +2c p\frac{\partial }{\partial p}r^{(12)}(z_1/z_2)
-c^*z\frac{\partial }{\partial z}r^{(12)}(z_1/z_2),
\\
&&c,c^*:\mbox{ central }, \quad [d,d^*]=0,\\
&&[d,\cL^\pm(z)]=-z\frac{\partial }{\partial z}\cL^\pm(z),
\\
&&[d^*,\cL^\pm(z)]=\pm 2 p\frac{\partial }{\partial p}\cL^+(z)
- z\frac{\partial }{\partial z}\widetilde{\cL}^-(z).
\end{eqnarray*}
Here we have set
\be
\widetilde{\cL}^-(z)&\!\!=\!\!&
h\otimes\left(\frac{1}{2}h+\sum_{n\neq 0}\frac{1}{1-p^n}z^n h^*_{-n}\right)
\nn\\
&&
+e\otimes\left(\sum_{n\in\Z}\frac{2}{1-w p^n}z^n f^*_{-n}\right)
+f\otimes\left(\sum_{n\in\Z}\frac{2}{1-w^{-1} p^n}z^n e^*_{-n}\right).
\en
Notice that, 
because of the quasi-Lie nature of $(\goth{g},\delta,\varphi)$, 
 $\goth{g}$ is a Lie subalgebra of the double $D\goth{g}$ 
whereas $\goth{g}^*$ is only a linear subspace 
(the Lie bracket does not close inside $\goth{g}^*$).

A similar description is possible for the classical limit of 
$\Aqp\bigl(\slth\bigr)$, but we omit the details.

\subsection{Comparison with Enriquez-Felder}\lb{sec:6.2}

In \cite{EF}, Enriquez and Felder studied the 
quasi-Hopf structure of an elliptic algebra $U_{\hbar}\goth{g}(\tau)$ 
associated with the face-type $R$ matrix. 
This algebra $U_{\hbar}\goth{g}(\tau)$ contains a central element $K$. 
Roughly speaking, $U_{\hbar}\goth{g}(\tau)$ and $\Uqp(\slth)$ 
are central extensions of the same algebra, 
as we already mentioned in section \ref{subsec:3.2}. 
Let us examine the main differences between the two algebras. 

The formulation of \cite{EF} starts with an 
elliptic curve with modulus $\tau$ and a coordinate $u$ on it.  
The latter plays the role of an `additive' spectral parameter, 
to be compared with our `multiplicative' spectral parameter $z=q^{2u}$. 
In the classical case, the relevant Manin pair in  \cite{EF}
is defined by assigning `positive' and `negative' 
parts {\it in powers of $u$}. 
Accordingly, in the construction of the half currents, 
the integration contours are chosen around a point in the $u$-plane. 
In our case, the integrations are taken along a circle around the origin 
on the $z$-plane. 

A more serious (perhaps related) difference arises in the quantum case. 
In \cite{EF}, the curve is fixed throughout quantization. 
In contrast, in the presence of the central term $c$, 
we have to deal simultaneously with 
two different elliptic curves with nomes $p$ and $\ps=pq^{-2c}$. 
This makes it difficult to adapt the geometric construction of \cite{EF}. 

The nature of the central element $K\in U_{\hbar}\goth{g}(\tau)$ 
and $c\in \Uqp(\slth)$ are quite different. 
Being dual to the grading element $d\in\g$, 
$K$ seems similar to the element $c^*\in\g^*$ in the double $D\g$ 
discussed in the previous subsection. 
For infinite dimensional representations and bosonization, 
we feel that it is more natural to consider extension by $c$.
A similar distinction has been discussed in the context of 
the double Yangian \cite{KLP3}.
In \cite{EF}, $U_\hbar\g(\tau)$ is initially endowed with 
a simple Hopf-algebra structure given by a Drinfeld-type coproduct. 
The quasi-Hopf structure related to the dynamical $RLL$ relation 
is then obtained by constructing a suitable twist from 
the initial Hopf structure. 
Such a Drinfeld-type coproduct persists in the presence of 
$c$ as well, but it provides only a quasi-Hopf structure 
(Appendix \ref{app:2}). 

Let us also mention a `physical' reason why we prefer $z$ to $u$. 
Recall Baxter's corner transfer matrix (CTM) method \cite{BaxBk}. 
A CTM is composed of a product of infinitely many Boltzmann weights. 
For the elliptic models, the individual weights 
(with appropriate normalization) are doubly periodic functions 
in $u$. 
The most important property of CTM is that, in the infinite lattice limit, 
the eigenvalues are all simple integral powers of $z$. 
This means that, in passing to the infinite lattice limit, 
one of the periods is lost. 
It is because the infinite lattice limit makes sense only inside 
the physical region, whereas 
shifting by another period takes us out of that region. 
Intuitively the $L$-operators generating the elliptic algebra 
are also products of Boltzmann weights on a single row of the lattice. 
The above argument indicates that in infinite dimensional representations 
the currents of the algebra possess only one period, 
and that $z$ is the natural variable to use. 

\subsection{Space of states}

Let us discuss how we view the space of states 
of the $k$-fusion unrestricted ABF model in connection with 
the algebras $\Bqla(\slth)$, $\Uqp(\slth)$. 
The parameter $r^*$ in $\lambda=(r^*+2)d+(s+1){1 \over 2} h$ 
corresponds to the elliptic nome. 
We shall argue below that the parameter $s$ corresponds to 
the boundary height degrees of freedom. 
 
First consider the `low temperature' limit $p,q\rightarrow 0$. 
Let us recall the `paths' of the spin $k/2$-XXZ model \cite{IIJMNT}.
A vertex-path $v$ is a semi-infinite sequence 
$v=(\cdots,{v(2)},{v(1)})$, where $v(l)\in  \{0,1,\cdots,k\}$. 
We have $k+1$ different ground state vertex-paths  
$\bar{v}_m$ ($m=0,1,\cdots,k$) given by 
\begin{equation}
  \bar{v}_m(l)=\biggl\{
  \begin{array}{ll}
    m&\mbox{for }l\equiv0 \mbox{ mod }2\\
    k-m&\mbox{for }l\equiv1 \mbox{ mod }2.
  \end{array}
\end{equation}
We say $v$ is an $m$-vertex-path if it satisfies the 
boundary condition $v(l)=\bar{v}_m(l)\quad \mbox{for } l\gg 1$. 
We assign a weight to an $m$-vertex-path by the formula 
\begin{eqnarray*}
{\rm wt}(v)=m(\Lambda_1-\Lambda_0)+ 
\sum_{l>0} (\bar{v}_m(l)-v(l))\alpha_1- h(v)(\alpha_0+\alpha_1), 
\end{eqnarray*}
where $h(v)$ is the `energy' of the path $v$ 
(see \cite{IIJMNT} for definition). 
The collection of all $m$-vertex-paths can be regarded as 
(the low temperature limit of) the space on which the CTM of the fusion 
vertex model acts. 
Its character $\sum_{v}q^{h(v)}z^{(wt(v),h)}$ (sum over the 
$m$-vertex-paths) is known to be the same as the character of the integrable 
highest weight $U_q(\slth)$-module $V(\mu_m)$ of highest weight 
\bea
\mu_m=(k-m)\Lambda_0+m\Lambda_1. 
\lb{mum}
\ena

Next we consider the paths for the unrestricted $k$-fusion ABF model. 
A face-path $s$ is a semi-infinite sequence $s=(\cdots,s(1),s(0))$ 
of integers $s(l)\in\Z$, 
subject to the admissibility condition 
$s(l)-s(l-1)\in\{k,k-2,\cdots,-k\}$ for $l\geq 1$. 
A face-path $s$ is called an $(m,n)$-face-path 
($m\in \{0,1,\cdots,k\}$, $n\in \Z$) if it satisfies the boundary condition 
$s(l)=\bar{s}_{m,n}(l)\quad \mbox{for } l\gg 1$, where 
$\bar{s}_{m,n}$ signifies the ground state face-path 
\begin{eqnarray*}
&&\bar{s}_{m,n}(0)=n + m,\\
&&\bar{s}_{m,n}(l)-\bar{s}_{m,n}(l-1)=2\bar{v}_m(l)-k. 
\end{eqnarray*}
{}From a face-path $s$, we can construct a vertex-path 
as $(\cdots, (s(2)-s(1)+k)/2,(s(1)-s(0)+k)/2)$, and 
conversely we obtain a unique face-path from a vertex-path 
up to a uniform shift $s(l)\rightarrow s(l)+a$ (for all $l$). 
Thus an $(m,n)$-face-path $s$ is uniquely represented 
by an $m$-vertex-path $v$ as  
\begin{eqnarray*}
&&{s}(0)=n + m+
\sum_{l>0} 2 (\bar{v}_m(l)-v(l)),\\
&&{s}(l)-{s}(l-1)=2 v(l)-k.
\end{eqnarray*}
The parameter $n$ determines the `boundary height at infinity', 
and the `bulk configuration' is described by some $m$-vertex-path $v$. 

Returning to the finite temperature situation $0<|p|<|q|<1$, let us consider 
the space ${\cal H}_{m,n}$ for the face model CTM 
under the boundary condition determined by $m,n$. 
As we have seen, for each fixed $n$, 
the character of ${\cal H}_{m,n}$ is the same as 
that of the $U_q(\slth)$-module $V(\mu_m)$. 
Since $\Bqla(\slth)$ has the same representations 
as the underlying algebra $U_q(\slth)$, 
${\cal H}_{m,n}$ can also be viewed as
the level $k$ irreducible highest weight module $V(\mu_m;r^*,s)$ 
over $\Bqla(\slth)$ with $\lambda=(r^*+2)d+(s+1){1 \over 2} h$.  
Let us consider the relation between $n$ and the parameter $s$. 
As we have discussed in the main text, 
the algebra $\Uqp(\slth)$ consists of two sectors, 
$U_q(\slth)$ and a Heisenberg algebra generated by $P,Q$ with $[P,Q]=-1$. 
A representation of this Heisenberg algebra 
is given by the zero mode lattice spanned by $e^{-nQ}|0\rangle$ ($n\in \Z$),
where $P|0\rangle=0$. 
The operator $P$ takes a fixed value $n$ on each state $e^{-nQ}|0\rangle$. 
Thus we can regard the direct sum 
\bea
\bigoplus_{n\in \Z\atop m=0,1,\cdots,k} V(\mu_m;r^*,n)\otimes e^{-nQ}
\lb{spst}
\ena
as a $\Uqp(\slth)$-module. 
This suggests that it is natural to identify $s$ with $n$. 
The above picture is consistent with the 
manner in which the dynamical shift appears
in the $RLL$ relation for $L^+(u,P)$. 

In summary, we identify the $\Uqp(\slth)$-module \eqref{spst} with 
the space ${\cal H}$,
\begin{equation}
  {\cal H}=\bigoplus_{n\in \Z\atop m=0,1,\cdots,k}{\cal H}_{m,n},
\end{equation}
on which the CTM of the $k$-fusion unrestricted ABF model acts 
under all possible boundary conditions. 

\subsection{Further issues}

Finally let us mention some related works and open problems. 

\begin{enumerate}
\item 
In \cite{Konno}, the algebra $\Uqp(\slth)$ was found by deforming the 
free field realization of the coset CFT 
$U({\slth})_l \otimes U({\slth})_k/U({\slth})_{l+k}$. 
As has been pointed out in \cite{LukPug2}, 
in the case $k=1$, $\Uqp(\slth)$ appears as the algebra of 
screening currents for the deformed 
Virasoro algebra (DVA).
(For the coset-type description of the RSOS model and DVA, 
see \cite{JMOh,JS}.) 
We wish to understand the conceptual meaning 
of DVA from 
the quasi-Hopf point of view. 
\item 
We note that the screening operators for the deformed $W_{n+1}$-
algebra coincide with  $E(u)$, $F(u)$ of $\Uqp(\slth)$
at level one up to some cocycle factors
which adjust signs in the commutation relations(Appendix A). 
A `higher rank extension' of DVA (the deformed $W_{n+1}$-algebra) 
and its screening currents have been studied in \cite{FeFr95,qWN}. 
The VO's for the $A^{(1)}_n$ face models 
is constructed in \cite{AJMP} by the use of the 
screening operator for the deformed $W_{n+1}$-algebra.
The cohomological structure of the Fock module of 
the deformed $W_{n+1}$-algebra is studied in \cite{FJMOP}. 
Even though not everything has yet been made clear, 
these works indicate that 
the deformed $W_{n+1}$-algebras play the role of the dynamical symmetry
of the $A^{(1)}_n$ face models. 
\item
In \cite{FrResh97}, the deformed $W$-algebra ${\cal W}_{q,t}(\bar\g)$ 
associated to an
arbitrary simple Lie algebra $\bar\g$ has been proposed. It
can be regarded as a quantization of the deformed Poisson $W$-algebra
obtained from a difference analogue of 
the Drinfeld-Sokolov reduction \cite{FRS,SS}.
On the other hand, we have obtained the algebra 
$U_{q,p}(\g)$ for an arbitrary non-twisted affine Lie algebra $\g$,
extending the results for
$U_{q,p}(\slth)$ as shown in Appendix A.
For non-simply laced $\g$, there is a considerable difference
between the Drinfeld-type currents for $U_{q,p}(\g)$ at level one and
the screening currents for ${\cal W}_{q,t}(\bar\g)$.
It seems natural to have such a difference since these
two have different CFT limits;
the former originates from the coset construction 
$U(\g)_k\otimes U(\g)_l/U(\g)_{k+l}$ and the
latter from the Drinfeld-Sokolov reduction of the 
loop group $G((z))$.
We thus expect the existence of another type of deformed $W$-algebra for 
non-simply laced $\goth{g}$ corresponding to $U_{q,p}(\g)$.
\item
A superalgebra version of the $W$-algebra was proposed recently 
in \cite{DF98}. We hope this class can be treated 
through the Lie superalgebra version of 
the quasi-Hopf structure \cite{ABRR}.
We note also that a construction of some extended version of DVA and 
deformed $W$-algebras from the 
vertex-type elliptic algebra $\Aqp(\widehat{\goth{sl}}_{n+1})$ 
is discussed in the works 
\cite{AFRS1,AFRS2,AFRS3,AFRS4}. 
\item
It is also tempting to guess that there is a 
`higher level extension' of the deformed Virasoro algebra 
whose screening currents are given by $E(u)$, $F(u)$ of $\Uqp(\slth)$
at level $k$. 
Such an algebra would be a deformation of 
(fractional) super-Virasoro algebra in CFT
(see \cite{Konno} for detailed discussions).
It is an open problem to find such an extended algebraic structure. 
\item 
In this paper as well as in \cite{JKOS1}, 
we focused attention to the case where the parameters $p,q$ are generic. 
In terms of lattice models, we considered only the unrestricted SOS case. 
The RSOS case corresponds to non-generic $p$ and needs 
special treatment. 
We wish to understand in particular the mechanism of 
obtaining possible extra singular vectors.\footnote{
In this connection we refer to the works \cite{BP97,BP98} 
on DVA and deformed $W$-algebras, 
where a detailed study is made on the 
Kac-determinant and properties of extra singular vectors at roots of unity.
}
The structure of the states of the RSOS model 
was also approached by using
the quantum affine algebra $U_q(\slth)$ \cite{JMOh}, and 
the DVA current was constructed within this language \cite{JMOh,JS}.
It is desirable to study the relationship
between this description and the one based on the quasi-Hopf algebra.
\item
An algebraic approach to the 
fusion ABF models has been presented on the basis of the 
quasi-Hopf algebra $\Bqla(\slth)$ and the elliptic algebra $\Uqp(\slth)$
in this work. 
Another interesting direction is 
to study Baxter's eight vertex model and Belavin's generalization.
Recently, a remarkable bosonization formula 
of the type I VO for the eight vertex model 
was proposed by Lashkevich and Pugai \cite{LP97}. 
They succeeded in reducing the problem to 
the already known bosonization for the ABF model 
through the use of intertwining vectors 
and Lukyanov's screening operators.  
To understand their bosonization scheme, 
it seems necessary to clarify 
the relationship between the intertwining vectors and the 
two twistors $F(\lambda)$ and $E(r)$, which 
define $\Bqla(\widehat{\goth{sl}}_n)$ and $\Aqp(\widehat{\goth{sl}}_n)$ 
respectively. 
It is also interesting to seek a more direct bosonization, 
which is intrinsically connected with the quasi-Hopf structure  
of $\Aqp(\slth)$ and does not rely on the bosonization of the ABF model. 
\end{enumerate}


\bigskip
\noindent
{\it Acknowledgments.}\quad 
We thank
Hidetoshi Awata, 
Jintai Ding, 
Benjamin Enriquez, 
Boris Feigin, 
Ian Grojnowski, 
Koji Hasegawa,
Hiroaki Kanno,
Harunobu Kubo, 
Michael Lashkevich, 
Tetsuji Miwa ,
Yaroslav Pugai, 
Takashi Takebe
and 
Jun Uchiyama
for discussions and interest.

\appendix
\setcounter{equation}{0}
\section{Elliptic currents for general $\gBig$}\lb{app:1}

In this appendix we give the elliptic currents and the algebra 
$\Uqp(\g)$ for non-twisted affine Lie algebras $\g$. 

\subsection{$U_q(\gbig)$}\lb{app:1.1}

Let $\g$ be an affine Lie algebra of non-twisted type associated with a  
generalized Cartan matrix $A=(a_{ij})$ and let $\cg$ be a corresponding 
simple finite dimensional Lie algebra.  Fixing an invariant bilinear form 
$(\ ,\ )$ on the Cartan subalgebra $\goth{h}$, we identify $\goth{h}^*$ with 
$\goth{h}$
via $(\ ,\ )$. Denoting the simple roots by $\alpha_j$, we set 
$b_{ij}=d_ia_{ij}=b_{ji}$ with $d_i=(\alpha_i,\alpha_i)/2$. Hence 
$B=(b_{ij})$ is  the symmetrized Cartan matrix. We also set $q_i=q^{d_i}$ and 
\begin{equation}
[n]_i=\frac{q_i^n-q_i^{-n}}{q_i-q_i^{-1}},\quad 
[n)_i=\frac{q^n-q^{-n}}{q_i-q_i^{-1}},\quad
\left[\matrix{m\cr n\cr}
\right]_i=\frac{[m]_i!}{[n]_i![m-n]_i!}.
\end{equation}

Consider the quantum affine algebra  $U_q(\g)$ of non-twisted type, and 
let 
\be
&&x_i^\pm(z)=\sum_{n\in \Z}x^\pm_{i,n} z^{-n},\\
&&\psi_i(q^{c/2}z)=q_i^{h_i}
\exp\left( (q_i-q_i^{-1}) \sum_{n>0} a_{i,n}z^{- n}\right),\\
&&
\varphi_i(q^{-c/2}z)=q_i^{-h_i}
\exp\left(-(q_i-q_i^{-1})\sum_{n>0} a_{i,-n}z^{ n}\right)
\en
be the Drinfeld currents ($i=1,2,..,{\rm rank}\ \cg$). 
The defining relations for $U_q(\g)$ read as follows: 
\be
&&c :\hbox{ central },\\
&& [h_i,d]=0,\quad [d,a_{i,n}]=n a_{i,n},\quad 
[d,x^{\pm}_{i,n}]=n x^{\pm}_{i,n}, \\
&&[h_i,a_{j,n}]=0,\qquad [h_i, x_j^\pm(z)]=\pm a_{ij} x_j^{\pm}(z),\\
&&
[a_{i,n},a_{j,m}]=\frac{[a_{ij}n]_i[c n)_j}{n}
q^{-c|n|}\delta_{n+m,0},\\
&&
[a_{i,n},x_j^+(z)]=\frac{[a_{ij}n]_i}{n}q^{-c|n|}z^n x_j^+(z),\\
&&
[a_{i,n},x_j^-(z)]=-\frac{[a_{ij}n]_i}{n} z^n x_j^-(z),
\\
&&(z-q^{\pm b_{ij}}w)
x_i^\pm(z)x_j^\pm(w)= (q^{\pm b_{ij}}z-w) x_i^\pm(w)x_j^\pm(z),
\\
&&[x_i^+(z),x_j^-(w)]=\frac{\delta_{i,j}}{q_i-q_i^{-1}}
\left(\delta\bigl(q^{-c}\frac{z}{w}\bigr)\psi_i(q^{c/2}w)
-\delta\bigl(q^{c}\frac{z}{w}\bigr)\varphi_j(q^{-c/2}w)
\right),\\
&&\sum_{\sigma\in S_a}\sum_{l=0}^a (-)^l
\left[\matrix{a\cr
l\cr}\right]_i
x^{\pm}_{i,m_{\sigma(1)}}\cdots x^{\pm}_{i,m_{\sigma(l)}}
x^{\pm}_{j,m} x^{\pm}_{i,m_{\sigma(l+1)}}\cdots x^{\pm}_{i,m_{\sigma(a)}}=0,
\nonumber\\
&&\qquad\quad (i\not=j,\ a=1-a_{ij},\ m_1,\cdots,m_a\in\Z).
\en
In the last line, $S_a$ denotes the symmetric group on $a$ letters.

\subsection{Elliptic currents}\lb{app:1.2}

Let us introduce the currents $u_i^\pm(z,p)\in U_q(\g)$ with $p=q^{2r}$ by 
\be
&&
u_i^+(z,p)=\exp\left(\sum_{n>0}\frac{1}{[r^*n)_i}a_{i,-n}(q^rz)^n\right),
\\
&&
u_i^-(z,p)=\exp\left(-\sum_{n>0}\frac{1}{[rn)_i}a_{i,n}(q^{-r}z)^{-n}\right).
\en
Then the following commutation relations hold. 
\begin{lem}\lb{lem:app1}
\begin{eqnarray}
&&u_i^+(z,p)x_j^+(w)=
\frac{(p^*q^{b_{ij}}z/w;p^*)_\infty}{(p^*q^{-b_{ij}}z/w;p^*)_\infty}
x_j^+(w)u_i^+(z,p),\\
&&u_i^+(z,p)x_j^-(w)=\frac{(p^*q^{-b_{ij}+c}z/w;p^*)_\infty}
{(p^*q^{b_{ij}+c}z/w;p^*)_\infty}
x_j^-(w)u_i^+(z,p),\\
&&u_i^-(z,p)x_j^+(w)=
\frac{(pq^{-b_{ij}-c}z/w;p)_\infty}{(pq^{b_{ij}-c}z/w;p)_\infty}
x_j^+(w)u_i^-(z,p),\\
&&u_i^-(z,p)x_j^-(w)=
\frac{(pq^{b_{ij}}z/w;p)_\infty}{(pq^{-b_{ij}}z/w;p)_\infty}
x_j^-(w)u_i^-(z,p),\\
&&u_i^+(z,p)u_j^-(w,p)=
\frac{(pq^{-c-b_{ij}}z/w;p)_\infty}{(pq^{-c+b_{ij}}z/w;p)_\infty}
\frac{(p^*q^{c+b_{ij}}z/w;p^*)_\infty}{(p^*q^{c-b_{ij}}z/w;p^*)_\infty}
u_j^-(w,p)u_i^+(z,p).
\end{eqnarray}
\end{lem}
Define the `dressed' currents $x_i^\pm(z,p)$, $\psi^{\pm}_{i}(z,p)$ in 
$U_q(\g)$ by 
\bea
&& x_i^+(z,p)=u_i^+(z,p)x_i^+(z), 
\lb{gdress1}
\\
&&x_i^-(z,p)=x_i^-(z)u_i^-(z,p),
\lb{gdress2}
\\
&&\psi^{+}_{i}(z,p)
=u_i^+(q^{c/2} z,p)\psi_i(z)u_i^-(q^{-c/2}z,p),
\lb{gdress3}
\\
&&\psi^{-}_{i}(z,p)
=u_i^+(q^{-c/2}z,p)\varphi_i(z)u_i^-(q^{c/2}z,p). 
\lb{gdress4}
\ena

Set $e_i(z)=x^+_i(z,p)$ and $f_i(z)=x^-_i(z,p)$. 
{}From  Lemma \ref{lem:app1}, we obtain 
\begin{prop} 
\bea
&&[h_i,a_{j,n}]=0,\quad [h_i,e_j(z)]=a_{ij} e_j(z),\quad 
[h_i,f_j(z)]=-a_{ij} f_j(z), 
\lb{B2}\\
&&[d,h_i]=0,\quad [d,a_{i,n}]= n a_{i,n},\quad \\
&&[d,e_i(z)]= -z\frac{\partial}{\partial z}e_i(z),\quad  
[d,f_i(z)]= -z\frac{\partial}{\partial z}f_i(z),
\lb{B3}\\
&&
[a_{i,n},a_{j,m}]=\frac{[a_{ij}n]_i[c n)_j}{n}
q^{-c|n|}\delta_{n+m,0},
\label{B4}\\
&&
[a_{i,n},e_j(z)]=\frac{[a_{ij}n]_i}{n}q^{-c|n|}z^n e_j(z),
\lb{B5}\\
&&
[a_{i,n},f_j(z)]=-\frac{[a_{ij}n]_i}{n} z^n f_j(z),
\lb{B6}\\
&& z\Theta_{p^*}\left(q^{b_{ij}} w/z\right) e_i(z)e_j(w)
=-w\Theta_{p^*}\left(q^{b_{ij}}z/w\right) 
e_j(w)e_i(z),
\lb{B7}\\
&& z \Theta_{p}\left(q^{-b_{ij}} w/z\right) f_i(z)f_j(w)
=-w \Theta_{p}\left(q^{-b_{ij}} z/w\right)
f_j(w)f_i(z),
\lb{B8}\\
&&[e_i(z),f_j(w)]=\frac{\delta_{i,j}}{q_i-q_i^{-1}}
\left(\delta\bigl(q^{-c}\frac{z}{w}\bigr)\psi^{+}_{i}(q^{c/2}w,p)
-\delta\bigl(q^{c}\frac{z}{w}\bigr)\psi^{-}_{i}(q^{-c/2}w,p)
\right),\lb{B9}\\
&&\sum_{\sigma\in S_a}\prod_{1\leq k<m\leq a}\frac{(p^*q^2 z_{\sigma(m)}/
z_{\sigma(k)};p^*)_{\infty}}{(p^*q^{-2} z_{\sigma(m)}/
z_{\sigma(k)};p^*)_{\infty}}\nonumber\\
&&\qquad\quad\times
\sum_{l=0}^a (-)^l
\left[\matrix{a\cr
l\cr}\right]_i \prod_{k=1}^l\frac{(p^*q^{b_{ij}} z/
z_{\sigma(k)};p^*)_{\infty}(p^*q^{-b_{ij}} 
z_{\sigma(k)}/z;p^*)_{\infty}}{(p^*q^{-b_{ij}} z/
z_{\sigma(k)};p^*)_{\infty}(p^*q^{b_{ij}} 
z_{\sigma(k)}/z;p^*)_{\infty}}\\
&&\qquad\quad\times e_{i}(z_{\sigma(1)})\cdots e_{i}(z_{\sigma(l)})
e_{j}(z) e_{i}(z_{\sigma(l+1)})\cdots e_{i}(z_{\sigma(a)})=0\qquad 
(i\not=j,\ a=1-a_{ij}),\nonumber\\
&&\sum_{\sigma\in S_a}\prod_{1\leq k<m\leq a}
\frac{(pq^{-2} z_{\sigma(k)}/
z_{\sigma(m)};p)_{\infty}}{(pq^{2} z_{\sigma(k)}/
z_{\sigma(m)};p)_{\infty}}\nonumber\\
&&\qquad\quad\times\sum_{l=0}^a (-)^l
\left[\matrix{a\cr
l\cr}\right]_i \prod_{k=1}^l\frac{(pq^{b_{ij}} z/
z_{\sigma(k)};p)_{\infty}(pq^{-b_{ij}} 
z_{\sigma(k)}/z;p)_{\infty}}{(pq^{-b_{ij}} z/
z_{\sigma(k)};p)_{\infty}(pq^{b_{ij}} 
z_{\sigma(k)}/z;p)_{\infty}}\\
&&\qquad\quad\times f_{i}(z_{\sigma(1)})\cdots f_{i}(z_{\sigma(l)})
f_{j}(z) f_{i}(z_{\sigma(l+1)})\cdots f_{i}(z_{\sigma(a)})=0\qquad 
(i\not=j,\ a=1-a_{ij}).\nonumber
\end{eqnarray}
\end{prop}

\subsection{$U_{q,p}(\gbig)$}\lb{app:1.3}

Let us introduce further a set of generators of the Heisenberg algebra 
$\{P_i,Q_i\}\ (i=1,2,..,{\rm rank }\ \cg)$ which 
commute with $U_q(\g)$ and satisfy 
\bea
&&[P_i, e^{Q_j}]=-\frac{a_{ij}}{2} e^{Q_j}. 
\ena
Setting $\tilde P_i=d_i P_i$,  $\tilde h_i=d_i h_i$ and  
\be
\!\!\!&&E_i(u)=e_i(z)e^{2Q_i}z^{-\frac{\tilde P_i-1}{r^*}},\\
\!\!\!&&F_i(u)=f_i(z) z ^{\frac{\tilde P_i+\tilde h_i-1}{r}},\\
\!\!\!&&H^{\pm}_i(z)=\psi^{\pm}_{i}(z)e^{2Q_i} 
(q^{\pm(r-\frac{c}{2})}z)^{-\frac{r-r^*}{rr^*}\tilde P_i
+\frac{1}{r}\tilde h_i},
\\
\!\!\!&&\hat{d}=d-\Delta^*+\Delta,\\
\!\!\!&&\Delta^*=
\frac{1}{2r^*}\sum_{i,j}(B^{-1})_{ij}(\tilde P_i-1)(\tilde P_j-3),
\quad \Delta=\frac{1}{2r}\sum_{i,j}(B^{-1})_{ij}(\tilde P_i+\tilde h_i-1)
(\tilde P_j+\tilde h_j-3),
\en
with 
$z=q^{2u}$, we have 
\bea
&&c:\hbox{ {\rm central}}, 
\lb{u1}
\\
&&[h_i,a_{j,n}]=0,\quad [h_i,E_j(u)]=a_{ij} E_j(u),\quad 
[h_i,F_j(u)]=-a_{ij} F_j(u), 
\lb{u2}\\
&&[\hat{d},h_i]=0,\quad [\hat{d},a_{i,n}]= n a_{i,n},\quad\\ 
&&[\hat{d},E_{i}(u)]=\left(-z\frac{\partial}{\partial z}+\frac{1}{r^*}\right)
E_i(u) , \quad
[\hat{d},F_{i}(u)]=\left(-z\frac{\partial}{\partial z}+\frac{1}{r}\right)
F_i(u) ,
\lb{u3}\\
&&
[a_{i,n},a_{j,m}]=\frac{[a_{ij}n]_i[c n)_j}{n}
q^{-c|n|}\delta_{n+m,0},
\label{u4}\\
&&
[a_{i,n},E_j(u)]=\frac{[a_{ij}n]_i}{n}q^{-c|n|}z^n E_j(u),
\lb{u5}\\
&&
[a_{i,n},F_j(u)]=-\frac{[a_{ij}n]_i}{n} z^n F_j(u),
\lb{u6}\\
&& \theta^*\left(u-v-\frac{b_{ij}}{2}\right) E_i(u)E_j(v)
=\theta^*\left(u-v+\frac{b_{ij}}{2}\right) 
E_j(v)E_i(u),
\lb{u7}\\
&& \theta\left(u-v+\frac{b_{ij}}{2} \right) F_i(u)F_j(v)
=\theta\left(u-v-\frac{b_{ij}}{2}\right)
F_j(v)F_i(u),
\lb{u8}\\
&&[E_i(u),F_j(v)]=\frac{\delta_{i,j}}{q_i-q_i^{-1}}
\left(\delta\bigl(q^{-c}\frac{z}{w}\bigr)H^+_{i}(q^{c/2}w)
-\delta\bigl(q^{c}\frac{z}{w}\bigr)H^-_{i}(q^{-c/2}w)
\right),\lb{u9}\\
&&\sum_{\sigma\in S_a}\prod_{1\leq k<m\leq a}z_{\sigma(k)}^{-\frac{2}{r^*}}
\frac{(p^*q^2 z_{\sigma(m)}/
z_{\sigma(k)};p^*)_{\infty}}{(p^*q^{-2} z_{\sigma(m)}/
z_{\sigma(k)};p^*)_{\infty}}\nonumber\\
&&\qquad\times
\sum_{l=0}^a (-)^l
\left[\matrix{a\cr
l\cr}\right]_i \prod_{k=1}^l \left(\frac{z}{z_{\sigma(k)}}
\right)^{\frac{b_{ij}}{r^*}}
\frac{(p^*q^{b_{ij}} z/
z_{\sigma(k)};p^*)_{\infty}(p^*q^{-b_{ij}} 
z_{\sigma(k)}/z;p^*)_{\infty}}{(p^*q^{-b_{ij}} z/
z_{\sigma(k)};p^*)_{\infty}(p^*q^{b_{ij}} 
z_{\sigma(k)}/z;p^*)_{\infty}}\\
&&\qquad\times E_{i}(u_{\sigma(1)})\cdots E_{i}(u_{\sigma(l)})
E_{j}(u) E_{i}(u_{\sigma(l+1)})\cdots E_{i}(u_{\sigma(a)})=0\qquad 
(i\not=j,\ a=1-a_{ij}),\nonumber\\
&&\sum_{\sigma\in S_a}\prod_{1\leq k<m\leq a}z_{\sigma(k)}^{\frac{2}{r}}
\frac{(pq^{-2} z_{\sigma(k)}/
z_{\sigma(m)};p)_{\infty}}{(pq^{2} z_{\sigma(k)}/
z_{\sigma(m)};p)_{\infty}}\nonumber\\
&&\qquad\times\sum_{l=0}^a (-)^l
\left[\matrix{a\cr
l\cr}\right]_i \prod_{k=1}^l\left(\frac{z}{z_{\sigma(k)}}
\right)^{\frac{b_{ij}}{r}}
\frac{(pq^{b_{ij}} z/
z_{\sigma(k)};p)_{\infty}(pq^{-b_{ij}} 
z_{\sigma(k)}/z;p)_{\infty}}{(pq^{-b_{ij}} z/
z_{\sigma(k)};p)_{\infty}(pq^{b_{ij}} 
z_{\sigma(k)}/z;p)_{\infty}}\\
&&\qquad\times F_{i}(u_{\sigma(1)})\cdots F_{i}(u_{\sigma(l)})
F_{j}(u) F_{i}(u_{\sigma(l+1)})\cdots F_{i}(u_{\sigma(a)})=0\qquad 
(i\not=j,\ a=1-a_{ij}).\nonumber
\end{eqnarray}

These are the generalizations of the defining relations of $\Uqp(\slth)$. 
Free field realizations of $U_{q,p}(\g)$ for
$\g = A^{(1)}_n,\ B^{(1)}_n$ and $ D^{(1)}_n$ are easily obtained. 
We will report on this subject in a future publication.

\medskip
\noindent
{\it Remark.}
Comparing  $U_{q,p}(\slnh)$ at level one 
with the relations among the screening currents of the deformed 
$W_{n}$-algebra\cite{FeFr95}, we have a difference by a sign factor 
$(-)^{a_{ij}}$. 
However, such 
discrepancy always occurs in the 
free field realization of (quantum) affine Lie algebras and is known to
be adjusted by using cocycle
factors\cite{FrKac} or by
using a procedure of central extension of the group algebra of the 
weight lattice. 
After such an adjustment, one can regard the algebra $U_{q,p}(\slnh)$ 
at level one as the algebra of the
screening currents of the deformed $W_{n}$-algebra.


\setcounter{equation}{0}
\section{Drinfeld coproduct}\label{app:2}

Besides the standard Hopf algebra structure, the quantum affine algebra 
$\uq$ is also endowed with the so-called Drinfeld coproduct: 
\bea
&&\Delta_\infty x=x\otimes 1+1\otimes x
\qquad (x=h,c,d),
\lb{copro0}\\
&&\Delta_\infty a_n= a_n\otimes 1+ q^{-c^{(1)}|n|}\otimes a_n, 
\lb{copro1}\\
&&\Delta_\infty x^+(z)=x^+(q^{-c^{(2)}}z)\otimes \psi^+(q^{-c^{(2)}/2}z)
+1 \otimes x^+(z),
\lb{copro2}\\
&&\Delta_\infty x^-(z)=x^-(z)\otimes 1 
+\psi^-(q^{-c^{(1)}/2}z)\otimes x^-(q^{-c^{(1)}}z). 
\lb{copro3}
\ena
The universal $R$ matrix ${\cal R}_\infty$ associated with this coproduct 
is given in \cite{DinKhoro}.
We have also an elliptic analog of the Drinfeld coproduct 
given as follows:
\bea
&&\Delta_{p,\infty} x=x\otimes 1+1\otimes x
\qquad (x=h,c,d),
\\
&&\Delta_{p,\infty} a_n=
\cases{
\displaystyle a_n\otimes 1+ \frac{[rn]}{[(r-c^{(1)})n]}\otimes a_n
& ($n>0$), \cr
\displaystyle \frac{[(r-c^{(1)}-c^{(2)})n]}
{[(r-c^{(1)})n]}a_n\otimes q^{c^{(2)}n}
+q^{c^{(1)}n} \otimes a_n
& ($n<0$), \cr}
\lb{Dco1}\\
&&\Delta_{p,\infty}
e(z,p)= e(q^{-c^{(2)}}z,p)
\otimes \psi^+(q^{-c^{(2)}/2}z,pq^{-2c^{(1)}})+1\otimes 
e(z,pq^{-2c^{(1)}}),  
\lb{Dco2}\\
&&\Delta_{p,\infty}
f(z,p)=f(z,p)\otimes 1+
\psi^-(q^{-c^{(1)}/2}z,p)\otimes f(q^{-c^{(1)}}z,pq^{-2c^{(1)}}). 
\lb{Dco3}
\ena
In terms of $\psi^\pm(z,p)$ we have 
\bea
&&
{\Delta}_{p,\infty}
\psi^\pm(z,p)= \psi^\pm(q^{\mp c^{(2)}/2}z,p)\otimes 
\psi^\pm(q^{\pm  c^{(1)}/2}z,pq^{-2c^{(1)}}).
\lb{Dco4}
\ena

As it turns out, this coproduct is obtained by a twist of 
\eqref{copro0}-\eqref{copro3}. 
Set 
\bea
F_\infty(p)=\exp\left(-\sum_{n>0}\frac{n}{[(r-c^{(1)})n][2n]}q^{rn}
a_{-n}\otimes a_n \right).
\lb{twist}
\ena
Then a simple computation shows the following. 
\begin{prop}
The twistor \eqref{twist} enjoys the shifted cocycle property 
\be
F_\infty^{(12)}(p)\left(\Delta_\infty\otimes \id\right)F_\infty(p)
=F_\infty^{(23)}(pq^{-2c^{(1)}})
\left(\id\otimes\Delta_\infty\right)F_\infty(p), 
\en
and satisfies 
\be
\Delta_{p,\infty}(a)=
F_\infty(p)\cdot \Delta_{\infty}(a)\cdot F_\infty(p)^{-1} 
\qquad \forall a\in \uq.
\en
\end{prop}
The universal $R$ matrix associated with this 
coproduct $\Delta_{p,\infty}$ is given by 
\bea
\cR_\infty(p)=F_\infty^{(21)}(p)\cR_\infty F_\infty(p)^{-1}.
\lb{Rinfr}
\ena
The classical limit of \eqref{Rinfr} reads 
\be
r_{\infty}(p)=
\frac{1}{2}h\otimes h+
2\sum_{n >0}\frac{p^n}{1-p^n}h_n\wedge h_{-n}
+2\sum_n e_n\otimes f_{-n} 
+c\otimes d+d\otimes c. 
\en
Upon skew-symmetrization, it gives the limiting case $w\rightarrow 0$ 
of the classical $r$ matrix \eqref{SOSr}. 
Note that the elliptic parameter $p$ enters only via the `Cartan' part. 
A similar classical $r$-matrix appeared also in the 
Drinfeld-Sokolov reduction \cite{FRS}. 


\setcounter{equation}{0}
\section{Evaluation modules}\label{app:3}

\subsection{Spin $l/2$ modules}\lb{app:3.1}

Let $l$ be a non-negative integer. 
We recall here the evaluation module of $\uq$ 
based on the spin $l/2$ representation. 

Let $V_l=\oplus_{m=0}^l \C v^l_m$, $V_{l,z}=V_l[z,z^{-1}]$. 
Define operators $h, S^\pm$ on $V_l$ by 
\bea
&&h v^l_m=(l-2m)v^l_m,
\qquad
S^\pm v^l_m=v^l_{m\mp1}, 
\lb{fd1}
\ena
where by convention we set $v^l_m=0$ for $m<0$ or $m>l$. 
In terms of the Drinfeld generators, 
the evaluation module $(\pi_{l,z},V_{l,z})$ is defined 
by the following formulas: 
\bea
&&\pi_{l,z}(c)=0, \qquad \pi_{l,z}(d)=z\frac{d}{dz}, 
\lb{fd4}\\
&&\pi_{l,z}(a_n)=
\frac{z^n}{n}\frac{1}{q-q^{-1}}
\left((q^n+q^{-n})q^{nh}-(q^{(l+1)n}+q^{-(l+1)n})\right),
\lb{fd5}\\
&&\pi_{l,z}(x^\pm(z'))
=S^\pm \Bigl[\frac{\pm h+l+2}{2}\Bigr] 
\delta\left(q^{h\pm 1}\frac{z}{z'}\right).
\lb{fd6}
\ena
In \eqref{fd6}, $[x]$ means $(q^x-q^{-x})/(q-q^{-1})$. 

The images of the elliptic currents \eqref{dress1}-\eqref{dress4}, 
\eqref{dress5} 
are then given as follows:
\bea
&&\pi_{l,z}(k(z'))=
\frac{\{q^{r-l}\frac{z}{z'}\}\{q^{r+l+2}\frac{z}{z'}\}}
{\{q^{r-l+2}\frac{z}{z'}\}\{q^{r+l+4}\frac{z}{z'}\}}
\frac{\{q^{r-l}\frac{z'}{z}\}\{q^{r+l+2}\frac{z'}{z}\}}
{\{q^{r-l+2}\frac{z'}{z}\}\{q^{r+l+4}\frac{z'}{z}\}}
\frac{(p;p)_\infty}{\Theta_p(q^{r-h}\frac{z'}{z})},
\lb{fd6'}\\
&&\pi_{l,z}(\psi^+(z',p))=
q^{h}
\frac{\Theta_p\left(q^{-l-1}\frac{z}{z'}\right)
\Theta_p\left(q^{l+1}\frac{z}{z'}\right)}
{\Theta_p\left(q^{-1+h}\frac{z}{z'}\right)
\Theta_p\left(q^{1+h}\frac{z}{z'}\right)},
\lb{fd7}\\
&&\pi_{l,z}(\psi^-(z',p))=
q^{-h}
\frac{\Theta_p\left(q^{-l-1}\frac{z'}{z}\right)
\Theta_p\left(q^{l+1}\frac{z'}{z}\right)}
{\Theta_p\left(q^{1-h}\frac{z'}{z}\right)
\Theta_p\left(q^{-1-h}\frac{z'}{z}\right)},
\lb{fd8}\\
&&
\pi_{l,z}(e(z',p))=S^+
\frac{q^{(- h-l)/2}}{1-q^{2}}
\frac{(q^{ h+l+2};p)_\infty(pq^{ h-l};p)_\infty}
{(p;p)_\infty(pq^{- 2};p)_\infty}
\delta\left(q^{h+ 1}\frac{z}{z'}\right),
\lb{fd91}\\
&&
\pi_{l,z}(f(z',p))=S^-
\frac{q^{(h-l)/2}}{1-q^{2}}
\frac{(q^{-h+l+2};p)_\infty(pq^{-h-l};p)_\infty}
{(p;p)_\infty(pq^{2};p)_\infty}
\delta\left(q^{h- 1}\frac{z}{z'}\right),
\lb{fd9}
\ena
where $\{z\}=(z;p,q^4)_\infty$.
In the text, we shall also write $V_{l,z}$ as $V_{l,u}$ with $z=q^{2u}$.

\subsection{$R$ matrix for spin $l/2$ representation} \lb{app:3.2}

Let $\cR^+=\cR$, $\cR^-={\cR^{(21)}}^{-1}$ be the universal $R$ matrices 
\footnote{Our $\cR$ here is ${\cR^{(21)}}^{-1}$ in \cite{JKOS1}, see
the remark at the end of section \ref{sec:2.1}.} of $\uq$, and let 
 $\cR^+(r,s)=\cR(r,s)$, $\cR^-(r,s)={\cR^{(21)}(r,s)}^{-1}$ be the elliptic 
counterparts. 
We consider their images in $V_{l,z_1}\otimes V_{m,z_2}$,
\be
R^{\pm}_{lm}(z_1/z_2)=\left(\pi_{l,z_1}
\otimes \pi_{m,z_2}\right)\cR^{\pm},
\qquad 
R^{\pm}_{lm}(z_1/z_2,s)
=\left(\pi_{l,z_1}\otimes \pi_{m,z_2}\right)\cR^{\pm}(r,s).
\en
The former has the form 
\be
R^{\pm}_{lm}(z)=\rho^{\pm}_{lm}(z)\bR_{lm}(z),
\en
where $\rho^\pm_{lm}(z)=\rho_{lm}(z^{\pm 1})^{\pm 1}$ is a scalar factor 
and $\bR_{lm}(z)$ is normalized as $\bR_{lm}(z)v^l_0\otimes v^m_0
=v^l_0\otimes v^m_0$. 
{}From the formula (3.10) of \cite{IIJMNT} for $\rho_{lm}(z)$, 
we find that 
\bea
&&R^{\pm}_{lm}(z,s)v^l_0\otimes v^m_0
=\rho^\pm_{lm}(z,p)v^l_0\otimes v^m_0, 
\nn\\
&&
\rho^+_{lm}(z,p)=
q^{lm/2}
\frac{\{pq^{l-m+2}z\}\{pq^{-l+m+2}z\}}{\{pq^{l+m+2}z\}\{pq^{-l-m+2}z\}}
\frac{\{q^{l+m+2}z^{-1}\}\{q^{-l-m+2}z^{-1}\}}
{\{q^{l-m+2}z^{-1}\}\{q^{-l+m+2}z^{-1}\}}.
\lb{rhop}
\ena
Set further $R^\pm_{lm}(z,s)=\rho^\pm_{lm}(z,p)\bR_{lm}(z,s)$. 
Noting that $\bR_{lm}(z)$ is a rational function in $z$, 
we find the following relation from (4.8) of \cite{JKOS1}:
\bea
\bR_{lm}(pz,s)=
q^{-\frac{1}{2}{h^{(1)}}^2-(s+h^{(2)})h^{(1)}}\cdot \bR_{lm}(z,s)
\cdot q^{\frac{1}{2}{h^{(1)}}^2+sh^{(1)}}.
\lb{pz}
\ena

Let us consider the image of the $L^+$ operator \eqref{Lplus} 
in the spin $l/2$ representation \eqref{fd6'}--\eqref{fd9}. 
With a suitable base change of the form $v^l_m\rightarrow g(h)v^l_m$, 
we find the following expression:
\bea
\pi_v(e^+(u,s))&\!\!\!=\!\!\!&
-S^+ \frac{\theta(\frac{l+h+2}{2})\theta(u-v-\frac{h+1}{2}-s)}
{\theta(u-v-\frac{h+1}{2})\theta(s)}
b(q^{h+1}w)^{-\frac{s}{r}}q^{-\frac{1}{r}\bigl(\frac{l+h+2}{2}\bigr)^2+
\frac{3}{4r}},
\lb{l/21}\\
\pi_v(f^+(u,s))&\!\!\!=\!\!\!&
S^- \frac{\theta(\frac{l-h+2}{2})\theta(u-v+\frac{h-1}{2}+s)}
{\theta(u-v-\frac{h-1}{2})\theta(s+h-1)}
b^{-1}(q^{h-1}w)^{\frac{s+h-1}{r}}
q^{-\frac{1}{r}\bigl(\frac{l-h+2}{2}\bigr)^2+\frac{1}{4r}},
\lb{l/22}\\
\pi_v(k^+_1(u,s))&\!\!\!=\!\!\!&
-\frac{\varphi_l(u-v-1)}{\theta(u-v-\frac{h+1}{2})}
w^{\frac{h}{2r}}q^{\frac{h^2-l(l+2)}{4r}},
\lb{l/23}\\
\pi_v(k^+_2(u,s))&\!\!\!=\!\!\!&
-\frac{\theta(u-v-\frac{h-1}{2})}{\varphi_l(u-v)}
w^{-\frac{h}{2r}}q^{-\frac{h^2-l(l+2)}{4r}}.
\lb{l/24}
\ena
Here $b$ is a constant, $w=q^{2v}$, and 
\bea
&&\varphi_l(u)=-z^{-\frac{l}{2r}}
\theta\bigl(u+\frac{l+1}{2}\bigr)\rho^+_{1l}(z,p)^{-1}.
\lb{varph}
\ena 
We note the relation 
\be
&&\varphi_l(u)\varphi_l(u-1)
=\theta(u-\frac{l+1}{2})\theta(u+\frac{l+1}{2}).
\en

We can get rid of 
the powers of $w$ and $q$ appearing in \eqref{l/21}-\eqref{l/24}
by the transformation \eqref{Lpr}.
Choosing 
\be
\mu(s,h)=w^{\frac{1}{2r}h(s+\frac{h}{2})}
q^{\frac{h^2-l(l+2)}{4r}s+\frac{h^3+(l+1)h^2-(l+1)^2h}{8r}}
\en
and $b=1$ we set 
\be 
R^+_{1l}(u-v,s)=(\id\otimes \pi_v)L^{+'}(u,s).
\en
The result is as follows:
\bea
R^+_{1l}(u,s)=\frac{1}{\varphi_l(u)}
\left(
\begin{array}{cc}
R_{++}(u) & R_{+-}(u) \\
R_{-+}(u) & R_{--}(u) \\
\end{array}
\right).
\lb{R1l1}
\ena
Here $\varphi_l(u)$ is given in \eqref{varph}, and 
the entries $R_{\ve,\ve'}(u)\in\End(V_{l,v})$ are given by 
\bea
&&
R_{++}(u) =-\frac{\theta(u+\frac{h+1}{2})
\theta(\frac{h-l}{2}+s)\theta(\frac{h+l}{2}+s+1)}
{\theta(s+h+1)\theta(s)},
\lb{R1l2}\\
&&
R_{+-}(u) =-S^-
\frac{\theta(\frac{l-h+2}{2})\theta(u+\frac{h-1}{2}+s)}{\theta(s+h-1)},
\lb{R1l3}\\
&&
R_{-+}(u) =S^+
\frac{\theta(\frac{l+h+2}{2})\theta(u-\frac{h+1}{2}-s)}{\theta(s)},
\lb{R1l4}\\
&&R_{--}(u)=-\theta(u-\frac{h-1}{2}).
\lb{R1l5}
\ena
In the simplest case $l=1$, this $R$-matrix coincides with \eqref{Rmat}
constructed from the image of twistors. 


\setcounter{equation}{0}
\section{Free field representation of $\Uqp(\slthBig)$}\lb{app:4}

In this section we review the free field representation of $\uq$. 
We then construct a free field representation of $U_{q,p}(\slth)$ 
following the prescription of section \ref{sec:3}. 

\subsection{Bosons}\lb{app:4.1}

Let $a_{n}$ be the bosons in section \ref{subsec:3.1}. 
In addition to them, we introduce two more kinds of bosons 
$a_{j,m}\ ( m\in\Z_{\not=0}\ j=1,2)$ satisfying 
the commutation relations
\bea
&&[a_{1,m},a_{1,n}]=\frac{[2m][(c+2)m]}{m}\delta_{m+n,0},\label{ccrb}\\
&&[a_{2,m},a_{2,n}]=-\frac{[2m][cm]}{m}\delta_{m+n,0}.\label{ccrc}
\ena
We  need also 
the zero-mode operators $Q_j$ and $P_j\ (j=0,1,2)$ satisfying
\bea
&&[P_0,Q_0]=-i,\qquad [P_1,Q_1]=2(c+2),\qquad [P_2,Q_2]=-2c.
\ena
It is also convenient to introduce the notation 
\begin{equation}
\al_+=\sqrt{\frac{r}{cr^*}},\qquad \al_-=-\sqrt{\frac{r^*}{cr}},\qquad
2\al_0=\al_+ + \al_-=\sqrt{\frac{c}{rr^*}}.
\end{equation}

Let us set 
\bea
h=-P_2, \qquad \al=\frac{1}{c}Q_2,\qquad \beta=-\sqrt{2}\al_-i Q_0.
\ena
Then $[h,\al]=2$.

We define the Fock space $\F_{J,\tilde M}$ by
\bea
&&\F_{J, \tilde M}=\bigoplus_{m,m'\in\Z}\F_{J,\tilde{M};m,m'},\\
&&\F_{J,\tilde{M};m,m'}=
\C[a_{-1},a_{-2},..;a_{j,-1},a_{j,-2},.. (j=1,2)]\nonumber\\
&&\qquad\qquad\qquad\otimes 
\C e^{\frac{J}{2(c+2)}Q_1}\otimes \C e^{m\beta}
\otimes \C e^{\frac{\tilde{M}}{2}\al+m'\al}.
\ena

\subsection{Bosonization of total currents}\lb{app:4.2}

Let us introduce the generating functions of bosons (boson fields),
\bea
{\phi}_j(A;B,C|z;D)&\!\!\!=\!\!\!&
-\frac{A}{BC}(Q_j+P_j\log z)+\tilde{\phi}_j(A;B,C|z;D),\\
\tilde{\phi}_j(A;B,C|z;D)&\!\!\!=\!\!\!&
\sum_{m\neq 0}\frac{[Am]}{[Bm][Cm]}a_{j,m}z^{-m}q^{D|m|},
\ena
and 
\begin{equation}
\phi_j^{(\pm)}(A;B|z;C)=\frac{P_j}{2}\log q+
(q-q^{-1})\sum_{m>0}\frac{[Am]}{[Bm]}a_{j,\pm m}z^{\mp m}
q^{C m}\quad (j=1,2).
\label{bosonfjpm}
\end{equation}
We sometimes use the abridgment
\bea
&&\phi_j(C|z;D)=\phi_j(A;A,C|z;D),\qquad \phi_j(C|z)=\phi_j(C|z;0). 
\ena

Now let us define the `parafermion fields' 
$\tilde{\Psi}(z)$ and $\tilde{\Psi}^{\dagger}(z)$ by $\tilde{\Psi}(z)=
\tilde{\Psi}^-(z),\ 
\tilde{\Psi}^{\dagger}(z)=\tilde{\Psi}^+(z) $, with 
\bea
&&\tilde{\Psi}^{\pm}(z)=\mp\frac{1}{(q-q^{-1})}\Bigl(
\tilde{\Psi}^\pm_I(z)-\tilde{\Psi}^\pm_{II}(z)
\Bigr),\label{pfcurrents}\\
&&\tilde{\Psi}^\pm_I(z)=:\EXP{\pm\tilde{\phi}_2(c|z;\pm\frac{c}{2})}
\EXP{-\phi^{(+)}_2\Bigl(1;2|z;\mp\frac{c+2}{2}\Bigr)\pm
\phi^{(+)}_1\Bigl(1;2|z;\mp\frac{c}{2}\Bigr)}:,\nonumber\\
&&\tilde{\Psi}^\pm_{II}(z)=
:\EXP{\pm\tilde{\phi}_2(c|z;\pm\frac{c}{2})}\EXP{\phi^{(-)}_2
\Bigl(1;2|z;\mp\frac{c+2}{2}\Bigr)\mp
\phi^{(-)}_1\Bigl(1;2|z;\mp\frac{c}{2}\Bigr)}:.\nonumber
\ena
Then we have 
\begin{prop}
The following currents $x^{\pm}(z)$ and operator $d$ 
with $h,\ c$ give a representation of $\uq$ on $\F_J=\F_{J,J}$:
\bea
&&x^+(z)= 
\tilde{\Psi}(z)\ :\EXP{-\sum_{n\neq 0}\frac{1}{[cn]}a_n z^{-n}}:e^\beta
\ e^\al,\label{uqe}\\
&&x^-(z)= 
\tilde{\Psi}^{\dagger}(z)\ 
:\EXP{\sum_{n\neq 0}\frac{q^{c|n|}}{[cn]}a_n z^{-n}}:e^{-\beta}\ e^{-\al} 
\label{uqf},\\
&&d=d_{1,2}+d_{a},\label{lzero}
\ena
where 
\bea
\!\!\!d_{1,2}&\!\!\!=\!\!\!&-\sum_{m>0}\frac{m^2}{[2m][(c+2)m]}a_{1,-m}a_{1,m}
+\sum_{m>0}\frac{m^2}{[2m][cm]}a_{2,-m}a_{2,m}-
\frac{P_1(P_1+2)}{4(c+2)},\\
\!\!\!d_{a}&\!\!\!=\!\!\!&-\sum_{m>0}\frac{m^2q^{cm}}{[2m][cm]}a_{-m}a_{m}.
\ena
\end{prop}
Note that this representation is slightly different from the one obtained 
by Matsuo\cite{Ma}. The main difference is in the identification of 
the Cartan operator $h$. See the discussion in subsection {D.3}. 

Note also that $\F_J$ gives a level $k$ highest weight representation of 
$U_q(\slth)$ for $c=k$ with the highest weight state
$$
|J\rangle=
1\otimes e^{{J\over 2(k+2)}Q_1}\otimes 1 \otimes 
e^{{J\over 2}\alpha}.
$$

With a substitution of \eqref{uqe} and \eqref{uqf} into \eqref{dress1} and 
\eqref{dress2},
the boson $a_n\ (n\in\Z_{\not=0})$, the currents $e(z,p), f(z,p)$ and 
$h,\ c,\ d$ give a representation of the elliptic currents of $\uq$
on $\F_J$. Explicitly, we have  
\begin{prop}
\bea
&&e(z,p)=
\tilde{\Psi}(z)\ :\EXP{-\sum_{n\neq 0}\frac{1}{[cn]}a_{0,n} z^{-n}}:e^\beta
\ e^\al,\label{bqpe}\\
&&f(z,p)= 
\tilde{\Psi}^{\dagger}(z)\ 
:\EXP{\sum_{n\neq 0}\frac{q^{c|n|}}{[cn]}a'_{0,n} z^{-n}}:e^{-\beta}\ e^{-\al} 
\label{bqpf}.
\ena
Here we introduced `dressed' bosons $a_{0,n}$ and $a'_{0,n}$ by
\bea
&&a_{0,n}=\cases{a_{n}& for $ n>0$\cr
                \displaystyle{\frac{[rn]}{[r^*n]}}q^{c|n|}
a_{n}&for $n<0$,
\cr}\\
&&a'_{0,n}=\frac{[r^*n]}{[rn]}a_{0,n}
\ena 
satisfying $[a_{0,m},a_{0,n}]=\frac{[2m][cm]}{m}\frac{[rm]}{[r^*m]}
\delta_{m+n,0}$ and 
$[a'_{0,m},a'_{0,n}]=\frac{[2m][cm]}{m}\frac{[r^*m]}{[rm]}\delta_{m+n,0}$.
\end{prop}

Let us next introduce the Heisenberg algebra generated by $P$ and $Q$. 
We realize them as
\begin{equation}
P-1=\sqrt{\frac{2rr^*}{c}}P_0+\frac{r^*}{c}h,\qquad
Q=-\sqrt{2}\alpha_0iQ_0.
\end{equation}
It is easy to check that $[Q,P]=1$ and that $P$ and $Q$ commute with $\uq$. 

Accordingly, we modify the Fock space $\F_J$ by $e^Q$ to $\hat\F_J$,
\begin{equation}
\hat{\F}_J=\bigoplus_{\mu\in\Z}\hat{\F}_{J,\mu},\qquad
\hat{\F}_{J,\mu}=\F_J\otimes e^{\mu Q}.
\end{equation}

Now we define the currents $K(z),\ E(z)$ and $F(z)$ by 
\eqref{U1}-\eqref{U3} replacing $e(z)$ and $f(z)$ with 
\eqref{bqpe} and \eqref{bqpf}, respectively.
Let us define also  
\bea
&&\hat{d}= d - \Delta_{-P+1,r^*} + \Delta_{-P-h+1,r}.
\ena
Then we have\cite{Konno} 
\begin{prop}
The currents  $K(z),\ E(z)$ and $F(z)$ and $h,\ c,\ \hat{d}$ give
a representation of $U_{q,p}(\slth)$ on $\hat{\F}_{J}$. Explicitly, 
these currents are given by 
\bea
&&K(z)=\,:\EXP{-\phi_0(1;2,r|z)}:,\label{currk}\\
&&E(z)= 
\Psi(z)\ :\EXP{-\phi_0(c|z)}:,\label{curre}\\
&&F(z)=
\Psi^{\dagger}(z)\ :\EXP{\phi'_0(c|z)}:,\label{currf}
\ena
where
\bea
&&\left.\matrix{\Psi^\dagger(z)\cr
\Psi(z)\cr}\right\}=\mp\frac{1}{(q-q^{-1})}\Bigl(
{\Psi}^\pm_I(z)-{\Psi}^\pm_{II}(z)
\Bigr),\label{fullpf}\\
&&{\Psi}^\pm_{I,II}(z)=\tilde{\Psi}^\pm_{I,II}(z)\ e^{\mp\al}\ 
z^{\pm\frac{1}{c}h},\nonumber
\ena
and
\bea
\phi_0(A;B,C|z;D)&\!\!\!=\!\!\!&\frac{A}{BC}\sqrt{\frac{2cr}{r^*}}
(i Q_0+P_0\log z)+\sum_{m\neq 0}\frac{[Am]}{[Bm][Cm]}
a_{0,m}z^{-m}
q^{D|m|},\label{bosonf0}\\
\phi'_0(A;B,C|z;D)&\!\!\!=\!\!\!&
\phi_0(A;B,C|z;D)\quad {\rm with }\ {r\leftrightarrow
 r^*,\  a_{0,m}\to a'_{0,m}}.
\ena
\end{prop}
Using the field $\phi_2(A|z;D)$,
the boson expression for the parafermion current $\Psi(z)$ 
(resp. $\Psi^\dagger(z)$) is obtained from the one for  $\tilde\Psi(z)$ 
(resp. $\tilde\Psi^\dagger(z)$) by
replacing the field $\tilde{\phi}_2(c|z;-c/2)$ with $\phi_2(c|z;-c/2)$
(resp. $\tilde{\phi}_2(c|z;c/2)$ with $\phi_2(c|z;c/2)$).

\Remark 
The parameterization of the vacuum charges of the Fock space
$\hat{\F}_{J,\tilde{M};m,m',\mu}=\F_{J,\tilde{M};m,m'}\otimes e^{\mu Q}$ 
is related to those of $\F_{J,M;n'n}$ in \cite{Konno}
as follows. Let us 
set $\al_{n',n}=\frac{1-n'}{2}\al_-+\frac{1-n}{2}\al_+$. Then
\begin{equation}
\tilde{M}+2m'=M,\qquad
m\beta+\mu Q=-\sqrt{2}\al_{n'n}iQ_0
\end{equation}
with $1-n'=2m+\mu$ and $1-n=\mu$. 
\medskip

\subsection{An alternative form}\lb{app:4.3} 

There is another way of constructing $\Uqp(\slth)$ from $\uq$
in terms of free bosons. 
Let us set 
\bea
&&
\bar{h}=-\sqrt{2c}P_0,\qquad \bar{\al}=-\sqrt{\frac{2}{c}}iQ_0.
\ena
Then $[\bar{h},\bar{\al}]=2$.

Define the Fock space $\bar{\F}_{J,\bar{M}}$ by
\bea
&&\bar{\F}_{J,\bar{M}}=
\bigoplus_{m,\bar{m}'\in\Z}\bar{\F}_{J,\bar{M};m,\bar{m}'},\\
&&\bar{\F}_{J,\bar{M};m,\bar{m}'}=
\C[a_{-1},a_{-2},..;a_{j,-1},a_{j,-2},.. (j=1,2)]\nonumber\\
&&\qquad\qquad\qquad\otimes
\C e^{\frac{J}{2(c+2)}Q_1}\otimes \C e^{\bar{m}'\al}
\otimes \C e^{\frac{\bar{M}}{2}\bar\al+m\bar{\al}}.
\ena

\begin{prop}
The following currents $x^{\pm}(z)$ and operator $\bar d$ 
with $\bar h,\ c$ give a representation of $\uq$ on 
$\bar\F_J=\bar{\F}_{J,J}$\cite{Ma}:
\bea
&&x^+(z)= 
{\Psi}(z)\ :\EXP{-\sum_{n\neq 0}\frac{1}{[cn]}a_n z^{-n}}:e^{\bar\al}
z^{\frac{1}{c}\bar h}
,\label{uqe2}\\
&&x^-(z)= 
{\Psi}^{\dagger}(z)\ 
:\EXP{\sum_{n\neq 0}\frac{q^{c|n|}}{[cn]}a_n z^{-n}}:e^{-\bar\al}
z^{-\frac{1}{c}\bar h} 
\label{uqf2},\\
&&\bar d=\bar d_{1,2}+\bar d_{a},\label{lzerobar}
\ena
where 
\bea
&&\bar d_{1,2}=-\sum_{m>0}\frac{m^2}{[2m][(c+2)m]}a_{1,-m}a_{1,m}
+\sum_{m>0}\frac{m^2}{[2m][cm]}a_{2,-m}a_{2,m}-
\frac{P_1(P_1+2)}{4(c+2)}+\frac{P_2^2}{4c},\nonumber\\
\\
&&\bar d_{a}=
-\sum_{m>0}\frac{m^2q^{cm}}{[2m][cm]}a_{-m}a_{m}-\frac{\bar{h}^2}{4c}.
\ena
\end{prop}
Then we have 
\begin{prop}
Dressing $x^\pm(z)$ by the procedure \eqref{dress1} and 
\eqref{dress2}, we have the following currents $e(z,p),f(z,p)$
with which the boson $a_n\ (n\in\Z_{\not=0})$ and 
$h,\ c,\ d$ give a representation of the elliptic currents of $\uq$ 
on $\bar{\F}_J$.  
\bea
&&e(z,p)=
{\Psi}(z)\ :\EXP{-\sum_{n\neq 0}\frac{1}{[cn]}a_{0,n} z^{-n}}:e^{\bar\al }
z^{\frac{1}{c}\bar h}
,\label{bqpe2}\\
&&f(z,p)= 
{\Psi}^{\dagger}(z)\ 
:\EXP{\sum_{n\neq 0}\frac{q^{c|n|}}{[cn]}a'_{0,n} z^{-n}}:e^{-\bar\al} 
 z^{-\frac{1}{c}\bar h}
\label{bqpf2}.
\ena
\end{prop}

In this case, we can obtain $\Uqp$ by dressing the elliptic currents via 
$\bar\al$ and $\bar h$
instead of adjoining $P$ and $Q$. This is a procedure of turning on the 
anomalous background charge $2\alpha_0$ in $\phi_0$. 
In conformal field theory, this corresponds to
the twist of the energy-momentum tensor by the Cartan operator. 
Then, the zero-mode lattice associated with $\bar{\al}$ gains one 
additional dimension and becomes 2-dimensional. Hence the Fock 
space $\bar\F_{J,\bar{M}}$ 
is changed to
\bea
&&\bar{\F}'_{J}=\bigoplus_{\bar{m}',\bar{n},\bar{n}'\in\Z}
\bar{\F}_{J; \bar{m}',\bar{n},\bar{n}'},\\
&&\bar{\F}'_{J; \bar{m}',\bar{n},\bar{n}'}=
\C[a_{-1},a_{-2},..;a_{j,-1},a_{j,-2},.. (j=1,2)]\nonumber\\
&&\qquad\qquad\qquad\otimes
\C e^{\frac{J}{2(c+2)}Q_1}\otimes \C e^{\bar{m}'\al}
\otimes \C e^{\left(\frac{\bar{n}}{2}\sqrt{\frac{r}{r^*}}+
\frac{\bar{n}'}{2}\sqrt{\frac{r^*}{r}}\right)\bar{\al}}.
\ena
 
\begin{prop}
The following currents  $\bar{K}(z)$, $\bar{E}(z)$, $\bar{F}(z)$ 
and $\tilde{d}$ with 
$h,\ c$ give
a representation of $U_{q,p}(\slth)$ on $\bar{\F}'_{J}$:
\bea
&&\bar K(z)=e^{\sqrt{c}\al_0\bar\al}\ k(z)\ z^{\frac{1}{2\sqrt{rr^*}}\bar h},
\label{currk2}\\
&&\bar E(z)=e^{-(1-\sqrt{\frac{r}{r^*}})\bar\al}\  e(z,p)\  z^{-\frac{1}{c}
(1-\sqrt{\frac{r}{r^*}})\bar h},\label{curre2}\\
&&\bar F(z)=e^{(1-\sqrt{\frac{r^*}{r}})\bar\al}\  f(z,p)\
z^{\frac{1}{c}(1-\sqrt{\frac{r^*}{r}})\bar h}\label{curr3},\\
&&\tilde{d}=\bar d-\frac{\al_0}{\sqrt{c}} \bar h .
\ena
\end{prop}

Expressing $P, \ \bar h $ by $P_j\ (j=0,1,2)$ and $\bar\al, \gamma$ 
by $Q_0,\ Q_2$,
the resultant $\bar K(z),\ \bar E(z),\ \bar F(z)$ and $\tilde{d}$ 
coincide with 
$K(z),\  E(z),\ F(z)$ and $\hat{d}$ in Proposition D.3, respectively.
The Fock space $\bar{\F}'_{J; \bar{m}',\bar{n},\bar{n}'}$ is isomorphic to
$\hat{\F}_{J,\tilde M;m,m',\mu}$ by 
\be
&&\bar{m}'=\frac{\tilde M}{2}+m',\quad \bar{n}=\mu,\quad \bar{n}'=-2m-\mu.
\en



\begin{thebibliography}{10}

\bibitem{DFJMN}
B.~Davies, O.~Foda, M.~Jimbo, T.~Miwa, and A.~Nakayashiki.
\newblock Diagonalization of the {XXZ} {Hamiltonian} by vertex operators.
\newblock {\em Comm. Math. Phys.}, 151:89--153, 1993.

\bibitem{JMMN}
M.~Jimbo, K.~Miki, T.~Miwa, and A.~Nakayashiki.
\newblock Correlation functions of the {X}{X}{Z} model for {${\Delta}<-1$}.
\newblock {\em Phys. Lett. A}, 168:256--263, 1992.

\bibitem{JM}
M.~Jimbo and T.~Miwa.
\newblock {\em Algebraic Analysis of Solvable Lattice Models}.
\newblock CBMS Regional Conference Series in Mathematics vol. 85, AMS, 1994.

\bibitem{FR}
I.~B. Frenkel and N.~Yu Reshetikhin.
\newblock Quantum affine algebras and holonomic difference equations.
\newblock {\em Comm. Math. Phys.}, 146:1--60, 1992.

\bibitem{BPZ}
A.~A. Belavin, A.~M. Polyakov, and A.~B. Zamolodchikov.
\newblock Infinite conformal symmetry in two-dimensional quantum field theory.
\newblock {\em Nucl. Phys.}, B241:333--380, 1984.

\bibitem{ABF}
G.~E. Andrews, R.~J. Baxter, and P.~J. Forrester.
\newblock Eight-vertex {SOS} model and generalized {Rogers-Ramanujan}-type
  identities.
\newblock {\em J. Stat. Phys.}, 35:193--266, 1984.

\bibitem{JMO}
M.~Jimbo, T.~Miwa, and M.~Okado.
\newblock Local state probabilities of solvable lattice models:an
  {$A^{(1)}_{n-1}$} family.
\newblock {\em Nucl. Phys.}, B300[FS22]:74--108, 1988.

\bibitem{DJKMO3}
E.~Date, M.~Jimbo, A.~Kuniba, T.~Miwa, and M.~Okado.
\newblock Exactly solvable {SOS} models {II}: Proof of the star-triangle
  relation and combinatorial identities.
\newblock {\em Adv. Stud. Pure Math.}, 16:17--122, 1988.

\bibitem{JMO3}
M.~Jimbo, T.~Miwa and M.~Okado.
\newblock Solvable lattice models related to the vector representation
of classical simple {Lie} algebras.
\newblock {\em Comm. Math. Phys.}, 116:507--525, 1988.

\bibitem{JMOh}
M.~Jimbo, T.~Miwa, and Y.~Ohta.
\newblock Structure of the space of states in {RSOS} models.
\newblock {\em Int. J. Mod. Phys.}, A8:1457--1477, 1993.

\bibitem{LukPug2}
S.~Lukyanov and Y.~Pugai.
\newblock Multi-point local height probabilities in the integrable {RSOS}
  model.
\newblock {\em Nucl. Phys.}, B473 [FS]:631--658, 1996.

\bibitem{MW96}
T.~Miwa and R.~Weston.
\newblock Boundary {ABF} models.
\newblock {\em Nucl. Phys.}, B 486[PM]:517--545, 1997.

\bibitem{qVir}
J.~Shiraishi, H.~Kubo, H.~Awata, and S.~Odake.
\newblock A quantum deformation of the {Virasoro} algebra and the {Macdonald}
  symmetric functions.
\newblock {\em Lett. Math. Phys.}, 38:33--57, 1996.

\bibitem{Konno}
H.~Konno.
\newblock An elliptic algebra {$U_{q,p}\bigl(\slth\bigr)$} and the fusion
  {RSOS} model, 1997.
\newblock {\em Comm. Math. Phys.}, 195:373-403, 1998.

\bibitem{FIJKMY}
O.~Foda, K.~Iohara, M.~Jimbo, R.~Kedem, T.~Miwa, and H.~Yan.
\newblock An elliptic quantum algebra for $\slth$.
\newblock {\em Lett. Math. Phys.}, 32:259--268, 1994.

\bibitem{Fel95}
G.~Felder.
\newblock Elliptic quantum groups.
\newblock {\em Proc. {ICMP Paris} 1994}, pages 211--218, 1995.

\bibitem{EF}
B.~Enriquez and G.~Felder.
\newblock Elliptic quantum groups {$E_{\tau,\eta}(\slth)$} and quasi-{Hopf}
  algebras, 1997.
\newblock q-alg/9703018.

\bibitem{Fron}
C.~Fr\o nsdal.
\newblock Generalization and exact deformations of quantum groups.
\newblock {\em Publ.RIMS, Kyoto Univ.}, 33:91--149, 1997.

\bibitem{Fron1}
C.~Fr\o nsdal.
\newblock Quasi-{Hopf} deformation of quantum groups.
\newblock {\em Lett. Math. Phys.}, 40:117--134, 1997.

\bibitem{QHA}
V.~G. Drinfeld.
\newblock Quasi-{Hopf} algebras.
\newblock {\em Leningrad Math. J.}, 1:1419--1457, 1990.

\bibitem{JKOS1}
M.~Jimbo, H.~Konno, S.~Odake, and J.~Shiraishi.
\newblock Quasi-{Hopf} twistors for elliptic quantum groups, 1997.
\newblock q-alg/9712029, to appear in Transformation Groups.

\bibitem{ABRR}
D.Arnaudon, E.Buffenoir, E.Ragoucy, and Ph.Roche.
\newblock Universal solutions of quantum dynamical {Yang-Baxter} equations,
1997.
\newblock q-alg/9712037.

\bibitem{DF}
J.~Ding and I.~Frenkel.
\newblock Isomorphism of two realizations of quantum affine algebra
  {$U_q\bigl(\widehat{\goth{gl}}_n\bigr)$}.
\newblock {\em Comm. Math. Phys.}, 156:277--300, 1993.

\bibitem{HouY}
B.~Hou and W.~Yang.
\newblock Dynamically twisted algebra {$A_{q,p;\hat{\pi}}\bigl(\slth\bigr)$} as
  current algebra generalizing screening currents of $q$-deformed {Virasoro}
  algebra, 1997.
\newblock q-alg/9709024.
%

\bibitem{Dri88}
V.~G. Drinfeld.
\newblock A new realization of {Yangians} and quantized affine algebras.
\newblock {\em Soviet Math. Dokl.}, 36:212--216, 1988.

\bibitem{DI97}
J.~Ding and K.~Iohara.
\newblock Generalization of {Drinfeld} quantum affine algebras.
\newblock {\em Lett. Math. Phys.}, 41:183--193, 1997.

\bibitem{KLP1}
S.~Khoroshkin, D.~Lebedev and S.~Pakuliak.
\newblock Elliptic algebra {${\cal A}_{q,p}(\slth)$} in the scaling limit.
\newblock {\em Comm. Math. Phys.}, 190:597--627, 1998.

\bibitem{KLP3}
S.~Khoroshkin, D.~Lebedev and S.~Pakuliak.
\newblock Yangian algebras and classical {Riemann} problems, 1997.
\newblock q-alg/9712057.

\bibitem{BaxBk}
R.~J. Baxter.
\newblock {\em Exactly Solved Models in Statistical Mechanics}.
\newblock Academic, London, 1982.

\bibitem{IIJMNT}
M.~Idzumi, K.~Iohara, M.~Jimbo, T.~Miwa, T.~Nakashima, and T.~Tokihiro.
\newblock Quantum affine symmetry in vertex models.
\newblock {\em Int. J. Mod. Phys.}, A8:1479--1511, 1993.

\bibitem{JS}
M.~Jimbo and J.~Shiraishi.
\newblock A coset-type construction for the deformed {Virasoro} algebra.
\newblock {\em Lett. Math. Phys.}, 43:173--185, 1998.

\bibitem{FeFr95}
B.~L. Feigin and E.~V. Frenkel.
\newblock Quantum {${\cal W}$}-algebras and elliptic algebras.
\newblock {\em Comm. Math. Phys.}, 178:653--678, 1996.

\bibitem{qWN}
H.~Awata, H.~Kubo, S.~Odake, and J.~Shiraishi.
\newblock Quantum {${\cal W}_N$} algebras and {Macdonald} polynomials.
\newblock {\em Comm. Math. Phys.}, 179:401--416, 1996.

\bibitem{AJMP}
Y.~Asai, M.Jimbo, T.Miwa, and Y.~Pugai.
\newblock Bosonization of vertex operators for the {$A^{(1)}_{n-1}$} face
model.
\newblock {\em J. Phys.} A29:6595--6616, 1996.

\bibitem{FJMOP}
B.~Feigin, M.~Jimbo, T.~Miwa, A.~Odesskii, and Y.~Pugai.
\newblock Algebra of screening operators for the deformed {$W_n$} algebra.
\newblock {\em Comm. Math. Phys.}, 191:501--541, 1998.

\bibitem{FrResh97}
E.~Frenkel and N.~Reshetikhin.
\newblock Deformations of ${\cal W}-$algebras associated to 
simple Lie algebras, 1997.
\newblock q-alg/9708006.

\bibitem{FRS}
E.~Frenkel, N.~Reshetikhin, M.A.~Semenov-Tian-Shansky.
\newblock Drinfeld-Sokolov reduction for difference 
operators and deformations of W-algebras I. The case of Virasoro
algebra.
\newblock {\em Comm. Math. Phys.}, 192:605--629, 1998.

\bibitem{SS}
M.A~ Semenov-Tian-Shansky and A.V.~Sevostyanov.
\newblock Drinfeld-Sokolov reduction for difference 
operators and deformations of W-algebras. II. General Semisimple Case.
\newblock {\em Comm. Math. Phys.}, 192:631--647, 1998.

\bibitem{DF98}
J.~Ding and B.Feigin.
\newblock Quantized {W}-algebra of {$\goth{sl}(2,1)$} : a construction from the
  quantization of screening operators, 1998.
\newblock math.QA/9801084.

\bibitem{AFRS1}
J.~Avan, L.~Frappat, M.~Rossi, and P.~Sorba.
\newblock Poisson structures on the center of the elliptic algebra 
{${\cal  A}_{q,p}({\hat sl}(2)_c)$}.
\newblock {\em Phys. Lett.} A235:323--334, 1997. 

\bibitem{AFRS2}
J.~Avan, L.~Frappat, M.~Rossi, and P.~Sorba.
\newblock New {$W_{q,p}\bigl(\goth{sl}(2)\bigr)$} algebras from the elliptic
  algebra {${\cal A}_{q,p}\bigl(\goth{sl}(2)_c\bigr)$}.
\newblock  {\em Phys. Lett.} A239:27--35, 1998. 

\bibitem{AFRS3}
J.~Avan, L.~Frappat, M.~Rossi, and P.~Sorba.
\newblock From quantum to elliptic algebras, 1997.
\newblock q-alg/9707034.

\bibitem{AFRS4}
J.~Avan, L.~Frappat, M.~Rossi, and P.~Sorba.
\newblock Deformed ${\cal W}_N$ algebras from elliptic $sl(N)$ algebras,
1998.
\newblock math.QA/9801105.

\bibitem{BP97}
P.~Bouwknegt and K.~Pilch.
\newblock The deformed {Virasoro} algebra at roots of unity, 1997.
\newblock q-alg/9710026.

\bibitem{BP98}
P.~Bouwknegt and K.~Pilch.
\newblock On deformed ${\cal W}$-algebras and quantum affine algebras,
1998.
\newblock math.QA/9801112.

\bibitem{LP97}
M.Lashkevich and Y.Pugai.
\newblock Free field construction for correlation functions of the eight-vertex
  model.
\newblock {\em Nucl. Phys.}, B516:623--651, 1998.

\bibitem{FrKac}
I.~Frenkel and V.G.~Kac.
\newblock Basic representations of affine Lie algebras and 
dual resonance models.
\newblock {\em Invent.Math.}, 62:23--66, 1980.

\bibitem{DinKhoro}
J.~Ding and S.~Khoroshkin.
\newblock Weyl group extension of quantized current algebras.
\newblock math.QA/9804139.

\bibitem{Ma}
A.~Matsuo.
\newblock A $q$-deformation of {Wakimoto} modules, primary fields and screening
  operators.
\newblock {\em Comm. Math. Phys.}, 161:33--48, 1994.

\end{thebibliography}

\end{document}